\documentclass[letterpaper, 10 pt, conference]{ieeeconf} 
\usepackage[]{times} 
\usepackage[T1]{fontenc}
\usepackage{graphicx,amstext,amsmath,amssymb,color,psfrag,mathabx,booktabs,mathtools,soul}
\usepackage{microtype}
\usepackage[utf8]{inputenc}
\usepackage[noadjust,nocompress]{cite}
\usepackage[colorlinks=true,    %no frame around URL
      urlcolor=black,     %no colors
      citecolor=blue,
      menucolor=black,    %no colors
      linkcolor=black,    %no colors 
      breaklinks=true,
]{hyperref}
\newif\iflong
\longtrue

\newcommand\mathieu[1]{{\color{magenta}{#1}}}
\newcommand{\Rl}[2]{\ensuremath{\mathbb{R}^{#1}_{#2}}}   % Euclidean space
\newcommand{\R}{\ensuremath{\mathbb{R}}}
\newcommand{\Rlp}{\ensuremath{\mathbb{R}_{>0}}}

\newcommand{\Rlo}{\ensuremath{\mathbb{R}_{\geq 0}}}

\newcommand{\Zo}{\ensuremath{\mathbb{Z}_{\geq 0}}}
\newcommand{\Zp}{\ensuremath{\mathbb{Z}_{> 0}}}
\newcommand{\Z}{\ensuremath{\mathbb{Z}}}

\definecolor{bleucit}{rgb}{0.2,0.4,0.6} % nouvelle couleur

\newcommand{\sabs}[1]{\ensuremath{\lceil {#1} \rfloor}}

\newcommand{\dst}{\displaystyle}

\definecolor{blue_cv}{rgb}{0.09,0.35,0.78}

\newcommand{\bm}[1]{\ensuremath{\mathbf{#1}}}

\newtheorem{defn}{Definition}

\newtheorem{prop}{Proposition}

\newtheorem{lem}{Lemma}

\newtheorem{thm}{Theorem}
\newtheorem{cor}{Corollary}

\newtheorem{rem}{Remark}
\newtheorem{sass}{Standing Assumption}
\usepackage{setspace}

\IEEEoverridecommandlockouts
\renewcommand\mathieu[1]{{{#1}}}
\begin{document}

\title{\LARGE Exploiting homogeneity for the optimal control of discrete-time systems: application to value iteration}

\author{Mathieu Granzotto, Romain Postoyan, Lucian Buşoniu, Dragan Nešić  and Jamal Daafouz \thanks{Mathieu Granzotto, Romain Postoyan and Jamal Daafouz are with the Universit\'e de Lorraine, CNRS, CRAN, F-54000 Nancy, France (e-mails: \{name.surname\}@univ-lorraine.fr).}\thanks{Lucian Bu\c{s}oniu is with the Department of Automation, Technical University of Cluj-Napoca, Memorandumului 28, 400114 Cluj-Napoca, Romania (e-mail: lucian.busoniu@aut.utcluj.ro).}\thanks{Dragan Nešić is with the Department of Electrical and Electronic Engineering, University of Melbourne, Parkville, VIC 3010, Australia (e-mail: dnesic@unimelb.edu.au). This work was supported by the Australian Research Council under the Discovery Project DP210102600 and by the France-Australia collaboration project IRP-ARS CNRS.}}
\date{}
\maketitle
%XXXX page numbers
%\thispagestyle{plain}
%\pagestyle{plain}
\IEEEpeerreviewmaketitle

\begin{abstract}
\mathieu{To investigate solutions of (near-)optimal control problems, we extend and exploit a notion of homogeneity recently proposed in the literature for discrete-time systems.} Assuming the plant dynamics \mathieu{is} homogeneous, we first derive a scaling property of \mathieu{its} solutions along rays provided the sequence of inputs is suitably modified. We then consider homogeneous cost functions and reveal how the optimal value function scales along rays. This result can be used to construct (near-)optimal inputs on the whole state space by only solving the original problem on a given compact manifold of a smaller dimension. {Compared to the related works of the literature,  we impose no conditions on the homogeneity degrees.} We demonstrate the strength of {this new result} by presenting a new approximate scheme for value iteration, which is one of the pillars of dynamic programming. \mathieu{The new} algorithm  provides guaranteed lower and upper estimates of the true value function at
any iteration and has several appealing features in terms of reduced computation. %First, it only requires to solve the original VI equation on a given compact manifold whose dimension is strictly smaller than the dimension of the state space. Second, the estimates of the value function are defined on the whole state space. These estimates are exact in some important, specific cases including minimum time and maximum hands-off control. 
A numerical case study is provided to illustrate the proposed algorithm.
\end{abstract}

%\noindent\romain{Comments are written in blue.}\\

%%%%%%%%%%%%%%%%%%%%%%%%%%%%%%%%%%%%%%%%
\vspace{-0.1em}\section{Introduction}\label{sect:intro}
%%%%%%%%%%%%%%%%%%%%%%%%%%%%%%%%%%%%%%%%
%More neutral, optimal control is .....
% it has known solutions in special cases
%it is hard for other type of costs  and nonlinear systems.
%in this paper we concentrate in discrete-time dynamical systems.... in this case....bellman next.
%because of the curse of dimensionality this is hard to do.
%when the real state spaces, a typical fashion is to employ compact state spaces. 

%say we only have approximations valued in the compact. we therefore need the set to be forward invariant, the design has no clear way to be defined, and is left for  the user.

%in this study we mitigate this effect by  homogeneity. With homogeneity
% we can reduce the study (we only have to solve the optimal control problem in a well defined compact set).

%while homogeneity in cont time has been explored for optimal control problems, refs, but there are much less work on discrete time (grimm et reihart). this is because the ct homogeneity def is not compatible (see article), which only works for degree =0. We use in this paper a notion introduced by %automatica sanchez et sanchez tac%, we need to extend this notion to system with inputs, we consider stage costs that are continuous, we derive the

%we exploit  theorem 1 may be used for DP and explicity MPC, and in fact we use this concepts in a VI scheme, and provide some guarantees. and see in 

%explicit MPC note

Optimal control seeks for   inputs   that   minimize   a given cost function along the \mathieu{trajectories} of the  considered plant dynamics, see, e.g., \cite{trelat2012optimal,Bertsekas-book12(adp),liberzon2011calculus}.  
While we know how to systematically construct optimal inputs for linear systems and quadratic costs since the 60's \cite{kalman1960contributions},  optimal control of general nonlinear systems \mathieu{with} general cost functions remains very challenging in general.
In this setting, dynamic programming  provides various algorithms to construct near-optimal inputs based on Bellman equation \cite{Bertsekas-book12(adp),bertsekas-tnnls15}.
When the state and input spaces are continuous, as in this work, dynamic programming algorithms cannot be implemented exactly. A typical procedure consists in considering compact subsets of the state and input spaces, respectively, which are then discretized to make the implementation of the algorithm tractable. A drawback of this approach is that we  obtain   \mbox{(near-)}optimal inputs for states  in the considered compact set \mathieu{only}. In addition, it is not always obvious that  a solution initialized in  this set would remain in it for all future times. Also, \mathieu{there is no clear recipe on how to select the compact set and typically this is done on a case by case basis by the user.}

%The price to pay is that this typically leads to  high computation cost\textcolor{black}{s}, which may be prohibitive. Indeed, because of the so-called curse of dimensionality, such approach is often intractable, in particular for problems with large dimensions and when the state and input spaces are continuous spaces.}   A typical \diff{way} to \diff{implement dynamic programming algorithms is to apply them on compact sets, see, e.g. \cite{Bertsekas-book12(adp),bertsekas-tnnls15,busoniu-et-al-aut10}}. However, \diff{by doing so, we only obtain (near-)optimal inputs for states in the considered compact set. It is then not always clear whether a solution in this set\textcolor{black}{, possibly the  optimal one,} would remain in it for all future times.} 
%Indeed, the optimal solution may not be confined to this compact set, and thus good approximations  are not guaranteed to be found. 
%One  then has to check if the chosen compact set is forward invariant, as solutions that step outside the set have consequently ill-defined value function approximations.   \diff{Also}, the \diff{choice} of the compact set is dependent of the problem and  has no definitive answer, being left to the user discretion.    %being -> hence?

In this work, we explain how these drawbacks can be addressed when the system and the cost functions satisfy \mathieu{homogeneity} properties. The underlying idea is that, in this case, we only have to solve the original optimal control problem on a given well-defined compact set  to derive \mbox{(near-)optimal} inputs for any state \mathieu{in the fully unbounded state space}. This approach was investigated in e.g.,  \cite{polyakov2020generalized,bhat2005geometric,coron2020model, nakamura2009homogeneous,mino1983homogeneity,rinehart-et-al-tac09(vi)} for continuous-time systems. The results on discrete-time systems we are aware of are more scarce \cite{Grimm-et-al-tac2005,rinehart-et-al-tac09(vi)} on the other hand, and only consider homogeneous properties of degree $1$ (or $0$ depending on the convention), which \mathieu{limits} their applicability. This restriction is due to the homogeneity notion used in \cite{Grimm-et-al-tac2005,rinehart-et-al-tac09(vi)}, which does not allow to derive scaling properties of the plant solutions along homogeneous rays when the degree is different from~$1$~{or~$0$}.  In this paper, we relax this degree constraint by exploiting a homogeneity definition recently proposed in the literature in \cite{Sanchez-et-al-ijrnc19,sanchez2020discrete}. In that way, we can derive scaling properties of the solutions for any homogeneity degree including those featured in  \cite[Proposition~2]{Grimm-et-al-tac2005} and \cite[proof of Corollary~1]{rinehart-et-al-tac09(vi)}  as {particular cases}. We are thus able to consider a much broader class of systems, which include polynomial systems for instance. %\st{we recover the results used in}

\mathieu{To this end, we} first extend the definition in \cite{Sanchez-et-al-ijrnc19,sanchez2020discrete} to be applicable for systems with inputs. We then derive a new scaling property of the system solutions, provided the sequence of inputs is suitably modified. Assuming the stage and terminal costs are also homogeneous, we explain how the optimal value function scales along homogeneous rays. As a result, by solving the original control problem at a given state, we obtain a solution to a \emph{different} optimal control problem by scaling the obtained sequence of inputs when moving along the corresponding ray. This differs from   \cite{polyakov2020generalized} where continuous-time systems are studied, as the results here employ scaled weighting constants, while   \cite{polyakov2020generalized} features instead   scaled cost horizons. \mathieu{Here,  we} deduce optimal sequences of inputs for the same cost function along the considered ray in some specific but important cases, which  include minimum-time control and maximum hands-off control \cite{nagahara-et-al-tac15,nagahara-et-al-max-hands-off-dt}, thereby covering and generalizing the results  of \cite{Grimm-et-al-tac2005}. 

%We believe that this result is important as it means that we only have to solve the original optimal control problem on a given compact set to derive optimal inputs for every point of the state space. 

We then explain how the general scaling property of the optimal value function can be exploited to obtain a new approximation scheme for value iteration (VI), which is one of the core algorithms of dynamic programming. This new scheme has the \mathieu{following} appealing features. First, it only requires to solve the original VI equation on a given compact manifold whose dimension is strictly smaller than the dimension of the state space. Second, the estimates of the value function are defined on the whole state space.  \mathieu{Finally, the} proposed algorithm is applied to the optimal control of a van der Pol oscillator.

The rest of the paper is organized as follows. The notation  is defined in  Section  \ref{sect:notation}.  The  main  results  are  given in  Section  \ref{sect-main-results}, where the homogeneous    scaling properties are stated.  The new approximation scheme for VI is presented in Section \ref{sect:dp}. The numerical case study is provided in  Section \ref{sect:ex}, and   Section \ref{sec:conc} concludes the paper. \iflong\else\mbox{ }The  proofs of Section \ref{sect:dp} and   \ref{sect-SA} have been  omitted for space reasons, and can be found in \cite{missing}.\fi

\vspace{-0.5em}
%%%%%%%%%%%%%%%%%%%%%%%%%%%%%%%%%%%%%%%%
\vspace{-0.1em}\section{Preliminaries}\label{sect:notation}
%%%%%%%%%%%%%%%%%%%%%%%%%%%%%%%%%%%%%%%%
Let $\R$ be the set of real numbers, $\Rlo:=[0,\infty)$, $\overline{\R}_{\geq 0}=[0,\infty]$, $\Rlp:=(0,\infty)$, $\Z$ be the set of integers, $\Zo:=\{0,1,2,\ldots\}$ and $\Zp:=\{1,2,\ldots\}$. Given a matrix $P\in\R^{n\times n}$ with $n\in\Zp$, we write $P\geq 0$ when $P$ is positive semi-definite  and we denote by $|x|_P:= x^\top P x$, for $x\in\R^{n}$, the $P$ induced norm. For any $d_1,\ldots,d_n\in\R$ with $n\in\Zp$, $\text{diag}(d_1,\ldots,d_n)$ stands for the $n\times n$ diagonal matrix, whose diagonal components are $d_1,\ldots,d_n$. The notation $(x,y)$ stands for $[x^{\top},\,y^{\top}]^{\top}$, where $x\in\Rl{n}{}$,  $y\in\Rl{m}{}$ and $n,m\in\Zp$. The Euclidean norm of a vector $x\in\Rl{n}{}$ with $n\in\Zp$ is denoted by $|x|$. We also use the so-called $\ell_0$ norm, which we denote $|x|_0$ for   $x\in\R^{n}$  with $n\in\Zp$, which is the number of components of $x$  that are not equal to $0$.  For any set $\mathcal{A}\subseteq\Rl{n}{}$, $x\in\R^{n}$ and $c\in\overline{\R}_{\geq 0}$, the map $\delta^{c}_{\mathcal{A}}:\R^{n}\to\overline{\R}_{\geq 0}$ is defined as $\delta^{c}_{\mathcal{A}}(x)=0$ when $x\in\mathcal{A}$ and $\delta^{c}_{\mathcal{A}}(x)=c$ when $x\notin\mathcal{A}$. Given vectors $x_1,\ldots,x_m\in\R^{n}$ with $n,m\in\Zp$, $\text{vect}(x_1,\ldots,x_m)$ is the vectorial space generated by the family $(x_1,\ldots,x_m)$. For any $x\in\R$ and $a\in\Rlp$, $\sabs{x}^{a}:=\text{sign}(x)|x|^{a}$. % A function $\chi:\Rlo\rightarrow\Rlo$ is of class $\mathcal{K}$ if it is continuous, zero at zero and
% strictly increasing, and it is of class $\Kinf$ if, in addition, it is unbounded. We say that a continuous function is of class $\mathcal{KL}$
% if for each $t\in\Rlo$, $\chi(\cdot,t)$ is of class $\mathcal{K}$, and, for each $s>0$, $\chi(s,\cdot)$ is decreasing to zero. 
We denote a sequence $\bm{u}_d=(u_0,\ldots,u_{d-1})$ of {(possibly infinitely many)} $d$ elements, with  $d\in\Zp\cup\{\infty\}$ and $u_i\in\R^{n}$ for any $i\in\{{0},\ldots,{d-1}\}$, $n\in\Zp$.   
% We use $\mathbb{I}$ to denote the identity matrix of appropriate dimension according to the context. 

The homogeneity definitions used in this work involve the next maps and sets.  Given $r=(r_1,\ldots,r_n)\in\Rlp^{n}$ with $n\in\Zp$, the \emph{dilation map} $\lambda^{r}_n:\Rlp\times\R^{n}\to\R^{n}$ is defined as $\lambda^r_n(\varepsilon,x)=(\varepsilon^{r_1}x_1,\ldots,\varepsilon^{r_n}x_n)$ for any  $\varepsilon\in\Rlp$ and $x=(x_1,\ldots,x_n)\in\R^{n}$. We denote $\lambda^r_n(\varepsilon,x)$ as $\lambda^{r}(\varepsilon)x$ when the dimension $n$ is clear from the context for the sake of convenience. We say that a dilation map is standard when $r=c(1,\ldots,1)\in\R^{n}$ for some $c\in\Rlp$. We call \emph{homogeneous ray} (induced by dilation map $\lambda^{r}$) passing at $x\in\R^{n}$, the set  $\mathcal{R}_r(x):=\{\lambda^{r}(\varepsilon)x\in\R^{n}\,:\,\varepsilon\in\Rlp\}$. When $x_1,x_2\in\mathcal{R}_r(x)$, we say that $x_1$ and $x_2$ belong to the same homogeneous ray. Given a set $\mathcal{A}\subseteq\R^{n}$ and $\varepsilon>0$, we also denote by $\lambda^{r}(\varepsilon)\mathcal{A}$ the set $\left\{y\in\R^{n}\,{:}\,\exists x\in\mathcal{A},  \quad y=\lambda^{r}(\varepsilon)x\right\}$. 
Given $r=(r_1,\ldots,r_n)\in\Rlp^{n}$ with $n\in\Zp$, $c\in\R$ and $k\in\Zp\cup\{\infty\}$, we denote by $\Lambda_{n,c,k}^{r}:\Rlp\times(\R^{n})^{k}\to(\R^{n})^{k}$  the map defined as $\Lambda_{n,c,k}^{r}(\varepsilon,\bm{u}_k):=\Big(\lambda^r(\varepsilon)u_0,\lambda^r(\varepsilon)^{c}u_1,\ldots,\lambda^r(\varepsilon)^{c^{k-1}}u_{k-1}\Big)$ for any $\bm{u}_k=(u_0,\ldots,u_{k-1})\in(\R^{n})^{k}$ and $\varepsilon\in\Rlp$. We write $\Lambda_{c,k}^{r}(\varepsilon)\bm{u}_k$ instead of $\Lambda_{n,c,k}^{r}(\varepsilon,\bm{u}_k)$ again when the dimension of $n$ is clear from context for the sake of convenience.

%%%%%%%%%%%%%%%%%%%%%%%%%%%%%%%%%%%%%%%%%%%%%%%%%
\vspace{-0.1em}\section{Main results}\label{sect-main-results}
%%%%%%%%%%%%%%%%%%%%%%%%%%%%%%%%%%%%%%%%%%%%%%%%%
In this section, we describe the homogeneity assumptions we make on the  system and the stage cost,
and how these allow to derive various scaling properties along homogeneous rays  of the solutions, the stage cost and finally of the  optimal value function. 
%%%%%%%%%%%%%%%%%%%%%%%%%%%%%%%%%
\vspace{-0.2em}\subsection{Homogeneous system}
%%%%%%%%%%%%%%%%%%%%%%%%%%%%%%%%%
Consider the discrete-time nonlinear system 
\begin{equation}
x^{+} = f(x,u),
\label{eq-plant}
\end{equation}
where $x\in\R^{n_x}$ is the plant state, $u\in\R^{n_u}$ is the control input and $n_x,n_u\in\Zp$. We assume throughout the paper that $f=(f_1,\ldots,f_{n_x}):\R^{n_x}\times\R^{n_u}\to\R^{n_x}$ satisfies the homogeneity property stated next.
%%%  
\begin{sass}[SA\ref{sass-homogeneity-f}]\label{sass-homogeneity-f} There exist $\nu\in\Rlp$, $r=(r_1,\ldots,r_{n_x})\in\Rlp^{n_x}$, and  $q=(q_1,\ldots,q_{n_u})\in\Rlp^{n_u}$ such that $f$ is homogeneous of degree $\nu\in\Rlp$ with respect to the dilation pair $(\lambda^r,\lambda^q)$ in the sense that for any $(x,u)\in\R^{n_x\times n_u}$ and $\varepsilon\in\Rlp$, $f_i(\lambda^{r}(\varepsilon)x,\lambda^{q}(\varepsilon)u)=\varepsilon^{r_i\nu}f_i(x,u)$, or, equivalently,  $f(\lambda^{r}(\varepsilon)x,\lambda^{q}(\varepsilon)u)=\lambda^{r}(\varepsilon)^{\nu}f(x,u)$. \hfill $\Box$
\end{sass}

The homogeneity property in SA\ref{sass-homogeneity-f} is inspired by  \cite[Definition 2]{Sanchez-et-al-ijrnc19} and extends it to vector fields with inputs. Some classes of vector fields $f$ satisfying SA\ref{sass-homogeneity-f} are provided in Section~\ref{subsect-sass-f}. Homogeneity definitions commonly found in the literature, as in, e.g.,  \cite{Grimm-et-al-tac2005,rinehart-et-al-tac09(vi)}, require $f_i(\lambda^{r}(\varepsilon)x,\lambda^{q}(\varepsilon)u)=\varepsilon^{r_i+\nu}f(x,u)$ instead of $f_i(\lambda^{r}(\varepsilon)x,\lambda^{q}(\varepsilon)u)=\varepsilon^{r_i\nu}f(x,u)$ as in SA\ref{sass-homogeneity-f}. The  homogeneity notions in  \cite[Definition 2]{Sanchez-et-al-ijrnc19} and SA\ref{sass-homogeneity-f} are better suited for vector fields arising in discrete-time dynamical systems. Indeed, they allow deriving scaling properties of the solutions as in  \cite[Lemma 3]{Sanchez-et-al-ijrnc19}, which is extended below to systems equipped with inputs. In terms of notation, we use $\phi(k,x,\bm{u}_{d})$ to denote the solution to system (\ref{eq-plant}) with initial condition $x$ and inputs sequence $\bm{u}_d\in{(\R^{n_u})^d}$ of length {$d\geq k$} at time $k\in\{0,\ldots,d\}$ {with the convention that $\phi(0,x,\cdot)=x$}. 

%{This  definition coincides at any time step $k\in\Zo$ to the  system's solution for a truncated input sequence, i.e., $\phi(k,x,\bm{u}_d)=\phi(k,x,\bm{u}_d|_k)$ where $\bm{u}_d|_k :=(u_0,\ldots,u_{k-1})$. While the latter is more explicit, the former keeps within reason the burden of notation.}

\begin{lem}\label{lem-scaling-solution} For any $x\in\R^{n_x}$, $\varepsilon>0$, $d\in\Zp\cup\{\infty\}$,  $\bm{u}_d\in{(\R^{n_u})^d}$ {and  $k\in\Zo$ with $k\leq d$,   }
\begin{equation}
\phi(k,\lambda^{r}(\varepsilon)x,\Lambda^{q}_{\nu,d}(\varepsilon)\bm{u}_{d}) =  \lambda^{r}(\varepsilon)^{\nu^{k}}\phi(k,x,\bm{u}_{d}).
\label{eq-lem-scaling}
\end{equation} 
\hfill $\Box$
\end{lem}

%%%%%%%%

\noindent\textbf{Proof:} Let $x\in\R^{n_x}$, $\varepsilon>0$ and $d\in\Zp\cup\{\infty\}$. {We proceed by induction. Considering the definition of $\phi$, we obtain
\begin{equation}
\phi(0,\lambda^{r}(\varepsilon)x,\Lambda^{q}_{\nu,{d}}(\varepsilon)\bm{u}_{d}) 
=  \lambda^{r}(\varepsilon)x,
\end{equation}
which corresponds to (\ref{eq-lem-scaling}) at $k=0$. Furthermore, by} the definitions of $\phi$ and $\Lambda^{q}_{\nu,{d}}${, it follows}  that
\begin{equation}
\begin{split}
\MoveEqLeft \phi(1,\lambda^{r}(\varepsilon)x,\Lambda^{q}_{\nu,{d}}(\varepsilon)\bm{u}_{d}) \\ 
&{}=  f\Big(\phi(0,\lambda^{r}(\varepsilon)x,\Lambda^{q}_{\nu,{d}}(\varepsilon)\bm{u}_{d}),\lambda^{q}(\varepsilon) u_0\Big)\\
&{}=  f\Big(\lambda^{r}(\varepsilon)x,\lambda^{q}(\varepsilon) u_0\Big).%\MoveEqLeft \phi(1,\lambda^{r}(\varepsilon)x,\Lambda^{q}_{\nu,{d}}(\varepsilon)\bm{u}_{d}) \\ &{}=   \phi(1,\lambda^{r}(\varepsilon)x,\lambda^{q}(\varepsilon)u_0) \\ 
%&{}=  f\Big(\phi(0,\lambda^{r}(\varepsilon)x,\lambda^{q}(\varepsilon)\bm{u}_1),\lambda^{q}(\varepsilon) u_0\Big)\\
%&{}=  f\Big(\lambda^{r}(\varepsilon)x,\lambda^{q}(\varepsilon) u_0\Big).
\end{split}
\end{equation}
We derive from SA\ref{sass-homogeneity-f}
\begin{equation}
\begin{split}
\phi(1,\lambda^{r}(\varepsilon)x,\Lambda^{q}_{\nu,{d}}(\varepsilon)\bm{u}_{d})  &{}=   \lambda^{r}(\varepsilon)^{\nu} f(x, u_0) \\ &{}= \lambda^{r}(\varepsilon)^{\nu}\phi(1,x,\bm{u}_{d}),
\end{split}    
\end{equation}
which corresponds to (\ref{eq-lem-scaling}) for $k=1$.

Suppose now that (\ref{eq-lem-scaling}) holds for ${k}\in\Zo$. We have
\begin{multline}
\phi(k+1,\lambda^{r}(\varepsilon)x,\Lambda^{q}_{\nu,d}(\varepsilon)\bm{u}_{d}) \\  {}= f\Big(\phi(k,\lambda^{r}(\varepsilon)x,\Lambda^{q}_{\nu,d}\bm{u}_{d}),\lambda^{q}(\varepsilon)^{\nu^{k}} u_{k}\Big).
\end{multline}
By assumption \mbox{$\phi(k,\!\lambda^{r}(\varepsilon)x,\!\Lambda^{q}_{\nu,d}(\varepsilon)\bm{u}_{d})\!=\!\lambda^{r}(\varepsilon)^{\nu^k}\!\phi(k,\!x,\!\bm{u}_{d})$}, thus
%\begin{equation}
\begin{multline}
\phi(k+1,\lambda^{r}(\varepsilon)x,\Lambda^{q}_{\nu,d}(\varepsilon)\bm{u}_{d})   \\ {}= f\Big(\lambda^{r}(\varepsilon)^{\nu^k}\phi(k,x,\bm{u}_{d}),\lambda^{q}(\varepsilon)^{\nu^k} u_{k}\Big).
\end{multline}
%\end{equation}
Noting that  $\lambda^{r}(\varepsilon)^{\nu^{k}}=\lambda^{r}(\varepsilon^{\nu^k})$ and $\lambda^{q}(\varepsilon)^{\nu^{k}}=\lambda^{q}(\varepsilon^{\nu^k})$ in view of the definition of $\lambda^{r}$,   we derive that by SA\ref{sass-homogeneity-f}
\begin{equation}
\begin{split}
\MoveEqLeft \phi(k+1,\lambda^{r}(\varepsilon)x,\Lambda^{q}_{\nu,d}(\varepsilon)\bm{u}_{d}) \\  &{}=f\Big(\lambda^{r}(\varepsilon^{\nu^k})\phi(k,x,\bm{u}_{d}),\lambda^{q}(\varepsilon^{\nu^k}) u_{k}\Big)\\
&{}= \lambda^{r}(\varepsilon^{\nu^{k}})^{\nu}f\Big(\phi(k,x,\bm{u}_{d}),u_{k}\Big)\\
&{}=
\lambda^{r}(\varepsilon)^{\nu^{k+1}}f\Big(\phi(k,x,\bm{u}_{d}),u_{k}\Big)\\
&{}= \lambda^{r}(\varepsilon)^{\nu^{k+1}}\phi(k+1,x,\bm{u}_{d}),
\end{split}
\end{equation}
which corresponds to the desired property, namely (\ref{eq-lem-scaling}) at $k+1$. We have proved that when (\ref{eq-lem-scaling}) holds for some $k\in\Zp$, it also does for $k+1$. The induction proof is complete.

\mbox{}\hfill $\blacksquare$ 
%%%%%%%%%

Lemma \ref{lem-scaling-solution} provides an explicit relationship between the solutions to system (\ref{eq-plant}) initialized on the same homogeneous ray (as defined in Section \ref{sect:notation}), provided their successive inputs are appropriately scaled by the map $\Lambda^{q}_{\nu,{d}}$.  

Next, we introduce the cost function and the  homogeneity properties we assume it satisfies.
%%% Jamal wondered bounded, answered.
%%%%%%%%%%%%%%%%%%%%%%%%%%%%%%%%%
\vspace{-0.2em}\subsection{Homogeneous cost}
%%%%%%%%%%%%%%%%%%%%%%%%%%%%%%%%%
We consider the cost function of horizon $d\in\Zp\cup\{\infty\}$, for any $x\in\R^{n_x}$ and any sequence of inputs of length $d$ $\bm{u}_{d}\in{(\R^{n_u})^d}$, 
\begin{multline}
J_{d,\gamma_1,\gamma_2}(x,\bm{u}_{d})  :=  \dst\sum_{k=0}^{d-1}\gamma_1^{\gamma_2^k}\ell\Big(\phi(k,x,\bm{u}_{d}),u_k\Big) \\ {}+ \gamma_1^{\gamma_2^d}\jmath(\phi(d,x,\bm{u}_{d})),
\label{eq-cost}
\end{multline}
where $\ell:\R^{n_x}\times\R^{n_u}\to\overline{\R}_{\geq 0}$ is the stage cost, $\jmath:\R^{n_x}\to\overline{\R}_{\geq 0}$ is the terminal cost, moreover   $\gamma_1\in\Rlp$ and $\gamma_2\in\R$  are weighting constants whose role is to capture the scaling effects due to homogeneity properties as it will be clarified in the sequel.  When $d=\infty$, $\jmath$ is the zero function.
The stage and terminal costs $\ell$ and $\jmath$ are assumed to satisfy the next homogeneity property.

\begin{sass}[SA\ref{sass-homogeneity-ell}]\label{sass-homogeneity-ell} There exists $\mu\in\R$ such that $\ell$ and $\jmath$ are homogeneous of degree $\mu$ with respect to dilation pairs $(\lambda^r,\lambda^q)$ and $\lambda^r$, respectively, in the sense that for any $(x,u)\in\R^{n_x\times n_u}$ and $\varepsilon\in\Rlp$, $\ell(\lambda^{r}(\varepsilon)x,\lambda^{q}(\varepsilon)u)=\varepsilon^{\mu}\ell(x,u)$ and $\jmath(\lambda^{r}(\varepsilon)x)=\varepsilon^{\mu}\jmath(x)$. \hfill $\Box$
\end{sass}

Examples of stage and terminal costs ensuring SA\ref{sass-homogeneity-ell} are given in Section \ref{subsect:sass-ell-jmath}. The next result is a consequence of Lemma~\ref{lem-scaling-solution} and SA\ref{sass-homogeneity-ell}.

\begin{lem}\label{lem-scaling-ell} For any $x\in\R^{n_x}$, $\varepsilon>0$, {$d\in\Zp$,  $\bm{u}_{d}\in{(\R^{n_u})^d}$ and $k\in\Zo$ with $k\leq d-1$, follows}
$\ell\!\left(\phi(k,\lambda^{r}(\varepsilon)x,\Lambda^{q}_{\nu,{d}}(\varepsilon)\bm{u}_{d},\lambda^{q}(\varepsilon)^{\nu^{k}}\!u_k\right)\! =  \varepsilon^{\mu\nu^k}\! \ell(\phi(k,x,\bm{u}_{d}),\allowbreak u_k)$ and
$\jmath\left(\phi({d},\lambda^{r}(\varepsilon)x,\Lambda^{q}_{\nu,{d}}(\varepsilon)\bm{u}_{d})\right) = \varepsilon^{\mu\nu^{d}}\jmath\left(\phi({d},x,\bm{u}_{d})\right)$
%\begin{equation}
%   \begin{split}
%       \MoveEqLeft\ell\left(\phi(k,\lambda^{r}(\varepsilon)x,\Lambda^{q}_{{{\nu,d}}(\varepsilon)\bm{u}|_{k}),\lambda^{q}(\varepsilon)^{\nu^{k}} u_k\right) \\ &{}=\varepsilon^{\mu\nu^k}\ell\left(\phi(k,x,\bm{u}|_{k}),u_k\right)\\
%\MoveEqLeft\jmath\left(\phi(k,\lambda^{r}(\varepsilon)x,\Lambda^{q}_{{{\nu,d}}(\varepsilon)\bm{u})\right) \\ &{}= %\varepsilon^{\mu\nu^k}\jmath\left(\phi(k,x,\bm{u})\right)\\
%   \end{split} 
%\label{eq:lem-scaling-ell}
%\end{equation}
where  $\nu$, $\mu$ respectively come from SA\ref{sass-homogeneity-f} and SA\ref{sass-homogeneity-ell}. \hfill $\Box$
\end{lem}

%%%%%
\iflong
\noindent\textbf{Proof:} Let $x\in\R^{n_x}$, $\varepsilon>0$, $d\in\Zp$, $k\in\{0,\ldots,d-1\}$ and $\bm{u}_{d}\in{(\R^{n_u})^d}$. In view of (\ref{eq-lem-scaling}), 
%\begin{multline}
$\ell\Big(\phi(k,\lambda^{r}(\varepsilon)x,\Lambda^{q}_{\nu,d}(\varepsilon)\bm{u}_{d}),\lambda^{q}(\varepsilon)^{\nu^{k}}u_{k}\Big) = \ell\Big(\lambda^{r}(\varepsilon)^{\nu^{k}}\phi(k,x,\bm{u}_{d}),\lambda^{q}(\varepsilon)^{\nu^k}u_{k}\Big)$.
%\end{multline}
Since $\lambda^{r}(\varepsilon)^{\nu^k}=\lambda^{r}(\varepsilon^{\nu^k})$ and $\lambda^{q}(\varepsilon)^{\nu^k}=\lambda^{q}(\varepsilon^{\nu^k})$, in view of SA\ref{sass-homogeneity-ell}, 
%\begin{multline}
$\ell\Big(\phi(k,\lambda^{r}(\varepsilon)x,\Lambda^{q}_{\nu,d}(\varepsilon)\bm{u}_{d}),\lambda^{q}(\varepsilon^{\nu^k})u_{k}\Big) = \left(\varepsilon^{\nu^k}\right)^{\mu}\ell\Big(\phi(k,x,\bm{u}_{d}),u_{k}\Big)$,
%\end{multline} 
which corresponds to the first equality in Lemma~\ref{lem-scaling-ell}. We similarly derive the second equality in Lemma~\ref{lem-scaling-ell}. \hfill $\blacksquare$
\fi
%%%%%%%%

We are ready to analyse the optimal value function associated to \eqref{eq-cost}.

\vspace{-0.2em}\subsection{Scaling property of the optimal value function}

Given $d\in\Zp\cup\{\infty\}$, $\gamma_1\in\Rlp$, $\gamma_2\in\R$ and $x\in\R^{n_x}$, we define the \emph{optimal value function} associated to cost $J_{d,\gamma_1,\gamma_2}$ at $x$ as
\begin{equation}
V_{d,\gamma_1,\gamma_2}(x) := \dst\min_{\bm{u}_{d}}J_{d,\gamma_1,\gamma_2}(x,\bm{u}_{d}),
\label{eq-V-gamma-d}
\end{equation}
which is assumed to exist, per the next assumption.

\begin{sass}[SA\ref{sass-existence-optimal-value-function}]\label{sass-existence-optimal-value-function} For any $d\in\Zp\cup\{\infty\}$, $\gamma_1\in\Rlp$, $\gamma_2\in\R$ and $x\in\R^{n_x}$, there exists a $d$-length sequence of inputs $\bm{u}_{{d}}^{\star}$ such that $V_{d,\gamma_1,\gamma_2}(x) = J_{d,\gamma_1,\gamma_2}(x,\bm{u}_{{d}}^{\star})<\infty$. \hfill $\Box$
\end{sass}

General conditions on $f$ and $\ell$ for  SA\ref{sass-existence-optimal-value-function} to hold can be found in \cite{Keerthi-Gilbert-tac85} for example.

The next theorem provides an explicit relationship on the optimal value function along   homogeneous rays. % \cite[Proposition 2]{Grimm-et-al-tac2005},  which requires $f$ to be of degree of homogeneity $1$, to the case where $f$ is of any degree of homogeneity as defined in SA\ref{sass-homogeneity-f}.
\begin{thm}\label{thm-scaling-optimal-value-function} For any $d\in\Zp\cup\{\infty\}$, $x\in\R^{n_x}$, denote $\bm{u}_{{d}}^{\star}$ an optimal sequence of inputs for the cost $J_{d,1,1}$ at $x$, then for any $\varepsilon>0$, \begin{equation}V_{d,\varepsilon^{-\mu},\nu}(\lambda^{r}(\varepsilon)x) = V_{d,1,1}(x),\label{eq-thm-scaling-optimal-value-function}\end{equation} %
where $\nu$ and $\mu$ respectively come from SA\ref{sass-homogeneity-f}-SA\ref{sass-homogeneity-ell}, and $\Lambda^{q}_{\nu,d}(\varepsilon)\bm{u}_{{d}}^{\star}$ is an optimal sequence of inputs for cost $J_{d,\varepsilon^{-\mu},\nu}$ at $\lambda^{r}(\varepsilon)x$, i.e., $J_{d,\varepsilon^{-\mu},\nu}(\lambda^{r}(\varepsilon)x,\Lambda^{q}_{\nu,d}(\varepsilon)\bm{u}_{{d}}^{\star})=V_{d,\varepsilon^{-\mu},\nu}(\lambda^{r}(\varepsilon)x)$, or, equivalently,
$V_{d,1,1}(\lambda^{r}(\varepsilon)x) = V_{d,\varepsilon^{\mu},\nu}(x)$, %\label{eq-thm-scaling-optimal-value-function-bis}
and, for $\bm{u}_{{d}}^{\star}$ an optimal sequence of inputs for costs $J_{d,1,1}$ at $\lambda^{r}(\varepsilon)(x)$,  $\Lambda^{q}_{\nu,d}(\varepsilon^{-1})\bm{u}_{{d}}^{\star}$ is an optimal sequence of inputs for cost $J_{d,\varepsilon^{\mu},\nu}$ at $x$, i.e., $J_{d,\varepsilon^{\mu},\nu}(x,\Lambda^{q}_{\nu,d}(\varepsilon^{-1})\bm{u}_{{d}}^{\star})=V_{d,\varepsilon^{\mu},\nu}(x)$. \mbox{}\hfill $\Box$ 
\end{thm}
% 

%%%%%%%%%%%%
 
\noindent\textbf{Proof:} Let $d\in\Zp\cup\{\infty\}$, $x\in\R^{n_x}$, $\bm{u}_{{d}}^{\star}$ be an optimal sequence of inputs for cost $J_{d,1,1}$ at $x$, which exists according to SA\ref{sass-existence-optimal-value-function}, and let $\varepsilon\in\Rlp$. From the definitions of $V_{d,1,1}$, $\lambda^{r}$, $\lambda^{q}$ and $\Lambda^q_{\nu,d}$, we derive that
\begin{align}
\MoveEqLeft V_{d,1,1}(x)=\dst\sum_{k=0}^{d-1}\ell\Big(\phi(k,x,\bm{u}_{{d}}^{\star}),u_k^{\star}\Big) +\jmath(\phi(d,x,\bm{u}_{{d}}^{\star})) \nonumber\\
&\hspace{-1.5em} {}=\dst\sum_{k=0}^{d-1}\ell\Big(\phi(k,\lambda^{r}(\varepsilon^{-1})\lambda^{r}(\varepsilon)x,\Lambda^{q}_{{\nu,d}}(\varepsilon^{-1})\Lambda^{q}_{{\nu,d}}(\varepsilon)\bm{u}_{{d}}^{\star}),\nonumber\\
&\hspace{-1.5em}{}\qquad\qquad\lambda^{q}(\varepsilon^{-1})^{\nu^k}\lambda^{q}(\varepsilon)^{\nu^k}u_k^{\star}\Big)\nonumber\\
& \hspace{-1em} {} + \jmath\left(\phi(d,\lambda^{r}(\varepsilon^{-1})\lambda^{r}(\varepsilon)x,\Lambda^{q}_{\nu,d}(\varepsilon^{-1})\Lambda^{q}_{\nu,d}(\varepsilon)\bm{u}_{{d}}^{\star})\right).
\end{align}
Hence, in view of Lemma \ref{lem-scaling-ell} and by definition of $V_{d,\varepsilon^{-\mu},\nu}(\lambda^{r}(\varepsilon)x)$, we have
\begin{align}
\MoveEqLeft V_{d,1,1}(x) \\ & \hspace{-2em} {}=  \dst\sum_{k=0}^{d-1}\varepsilon^{-\mu\nu^k}\ell\Big(\phi(k,\lambda^{r}(\varepsilon)x,\Lambda^{q}_{{\nu,d}}(\varepsilon)\bm{u}_{{d}}^{\star}),\lambda^{q}(\varepsilon)^{\nu^k}u_k^{\star}\Big)\nonumber\\
&\hspace{-2em} {}+ \varepsilon^{-\mu\nu^d} \jmath\left(\phi(d,\lambda^{r}(\varepsilon)x,\Lambda^{q}_{\nu,d}(\varepsilon)\bm{u}_{{d}}^{\star})\right) \geq  V_{d,\varepsilon^{-\mu},\nu}(\lambda^{r}(\varepsilon)x).\nonumber\label{eq-proof-thm-scaling-optimal-value-function-lower-bound}  
\end{align}
On the other hand, 
\begin{align}
V_{d,\varepsilon^{-\mu},\nu}(\lambda^{r}(\varepsilon)x) &= \dst\sum_{k=0}^{d-1}\varepsilon^{-\mu\nu^{k}}\ell\Big(\phi(k,\lambda^{r}(\varepsilon)x,\bm{v}_{{d}}^{\star}),v_k^{\star}\Big)\nonumber\\
 & \quad {}+ \varepsilon^{-\mu\nu^{d}}\jmath(\phi(d,\lambda^{r}(\varepsilon)x),\bm{v}_{{d}}^{\star}), 
\end{align}
where $\bm{v}_{{d}}^{\star}$ is an optimal sequence of inputs for cost $J_{d,\varepsilon^{-\mu},\nu}$ at initial state $\lambda^{r}(\varepsilon)x$, which exists by SA\ref{sass-existence-optimal-value-function}. We have that
\begin{align}
\MoveEqLeft V_{d,\varepsilon^{-\mu},\nu}(\lambda^{r}(\varepsilon)x)\nonumber\\ &{} =  \dst\sum_{k=0}^{d-1}\varepsilon^{-\mu\nu^{k}}\ell\Big(\phi(k,\lambda^{r}(\varepsilon)x,\Lambda^{q}_{{\nu,d}}(\varepsilon)\Lambda^{q}_{{\nu,d}}(\varepsilon^{-1})\bm{v}_{{d}}^{\star}),\nonumber\\
& \qquad\qquad\qquad\quad\lambda^{q}(\varepsilon)^{\nu^k}\lambda^{q}(\varepsilon^{-1})^{\nu^k}v_k^{\star}\Big)\\
 &\quad {} + \varepsilon^{-\mu\nu^{d}}\jmath\Big(\phi(d,\lambda^{r}(\varepsilon)x,\Lambda^{q}_{\nu,d}(\varepsilon)\Lambda^{q}_{\nu,d}(\varepsilon^{-1})\bm{v}_{{d}}^{\star})\Big),\nonumber
\end{align}
from which we derive, in view of Lemma \ref{lem-scaling-ell} and the definition of $V_{d,1,1}(x)$,
\begin{align}
\MoveEqLeft V_{d,\varepsilon^{-\mu},\nu}(\lambda^{r}(\varepsilon)x) \nonumber\\&= \dst\sum_{k=0}^{d-1}\ell\Big(\phi(k,x,\Lambda^{q}_{{\nu,d}}(\varepsilon^{-1})\bm{v}_{{d}}^{\star}),\lambda^{q}(\varepsilon^{-1})^{\nu^k}v_k^{\star}\Big)\nonumber\\
& \quad {}+ \jmath\Big(\phi(d,x,\Lambda^{q}_{\nu,d}(\varepsilon^{-1})\bm{v}_{{d}}^{\star}))\Big) \geq  V_{d,1,1}(x).
\label{eq-proof-thm-scaling-optimal-value-function-upper-bound}
\end{align}
We deduce from (\ref{eq-proof-thm-scaling-optimal-value-function-lower-bound}) and (\ref{eq-proof-thm-scaling-optimal-value-function-upper-bound}) that $V_{d,\varepsilon^{-\mu},\nu}(\lambda^{r}(\varepsilon)x) = V_{d,1,1}(x)$. As a consequence, in view of the equality in (\ref{eq-proof-thm-scaling-optimal-value-function-lower-bound}),  $\Lambda^{q}_{\nu,d}{(\varepsilon)}\bm{v}_{{d}}^{\star}$ is an optimal sequence of inputs for cost $J_{d,\varepsilon^{-\mu},\nu}$ at $\lambda^{r}(\varepsilon)x$. 

We now prove the last part of Theorem \ref{thm-scaling-optimal-value-function}. Let $z:=\lambda^{r}(\varepsilon)x$. In view of the first part, %(\ref{eq-thm-scaling-optimal-value-function}),
$V_{d,1,1}(z) = V_{d,\varepsilon^{\mu},\nu}(\lambda^{r}(\varepsilon^{-1})z)$,
which, by using the definition of $z$, is equivalent to
$V_{d,1,1}(\lambda^{r}(\varepsilon)x) = V_{d,\varepsilon^{\mu},\nu}(x)$.
We have obtained the desired result.%(\ref{eq-thm-scaling-optimal-value-function-bis}). 
\hfill $\blacksquare$ 
%%%%%%%%%%%

Theorem \ref{thm-scaling-optimal-value-function} provides a scaling-like property of the optimal value function along homogeneous rays. In particular, by computing $V_{d,1,1}$ and an associated sequence of optimal inputs at some $x\in\R^{n_x}\setminus\{0\}$, we can deduce a scaled version of the optimal value function and an associated sequence of optimal inputs for any point on the   homogeneous ray $\mathcal{R}_r(x)$. However, the obtained optimal value function differs in general in Theorem \ref{thm-scaling-optimal-value-function}, except in two specific cases   formalized in the next corollary,
the second one being novel, which is a direct application of Theorem~\ref{thm-scaling-optimal-value-function}.

\begin{cor}\label{cor-scaling-optimal-value-function} Suppose that either SA\ref{sass-homogeneity-f} holds with $\nu=1$ or SA\ref{sass-homogeneity-ell} holds with $\mu=0$. Then, for any $d\in\Zp\cup\{\infty\}$, $x\in\R^{n_x}$, denote $\bm{u}_{{d}}^{\star}$ an optimal sequence of inputs for cost $J_{d,1,1}$ at $x$, then for any $\varepsilon>0$ \begin{equation}
V_{d,1,1}(\lambda^{r}(\varepsilon)x)  =  \varepsilon^{\mu}V_{d,1,1}(x),
\label{eq-cor-scaling-optimal-value-function}
\end{equation}
and $\Lambda^{q}_{1,d}(\varepsilon)\bm{u}_{{d}}^{\star}$ is an optimal sequence of inputs for cost $J_{d,1,1}$ at $\lambda^{r}(\varepsilon)x$.\hfill $\Box$
\end{cor}
\noindent\textbf{Proof:} { Let $x\in\R^{n_x}$ and $\varepsilon>0$.
When SA\ref{sass-homogeneity-f} holds with $\nu=1$ and $\mu\in\R$, (\ref{eq-thm-scaling-optimal-value-function}) gives  
\begin{equation}
 V_{d,\varepsilon^{-\mu},1}(\lambda^{r}(\varepsilon)x)  =  \varepsilon^{-\mu}V_{d,1,1}(\lambda^{r}(\varepsilon)x).
\end{equation}
Hence, in view of (\ref{eq-thm-scaling-optimal-value-function}),
\begin{equation}
 V_{d,1,1}(\lambda^{r}(\varepsilon)x) = e^{\mu} V_{d,1,1}(x).
\end{equation}

 When SA\ref{sass-homogeneity-ell} holds with $\mu=0$ and $\nu\in\Rlp$, (\ref{eq-thm-scaling-optimal-value-function}) gives  
 \begin{equation}
 V_{d,1,\nu}(\lambda^{r}(\varepsilon)x)  = V_{d,1,1}(\lambda^{r}(\varepsilon)x).
 \end{equation}
 Hence, in view of (\ref{eq-thm-scaling-optimal-value-function}),
 \begin{equation}
V_{d,1,1}(\lambda^{r}(\varepsilon)x)  =  V_{d,1,1}(x)
 \end{equation}
 and the proof is complete.
% \romain{This result may be useful for minimum time optimal control. Indeed, in this  case, when the target is the origin, we have $\ell(x,u)=\bm{1}_{0}(x)$ where $\bm{1}_{0}$ is the indicator function of the origin. Then, $\ell$ is homogeneous of degree $\mu=0$ with respect to any dilation pair according to SA\ref{sass-homogeneity-ell}. 
}
\mbox{} \hfill $\blacksquare$ 

In Corollary \ref{cor-scaling-optimal-value-function}, the optimal value function is proportional along homogeneous rays under extra conditions compared to Theorem \ref{thm-scaling-optimal-value-function}. As a result, by computing the optimal value function on any $n_x$-sphere of $\R^{n_x}$ centered at the origin for instance and an associated sequence of optimal inputs, we {can} obtain the value of the optimal value function as well as a sequence of optimal inputs for any point in $\R^{n_x}\setminus\{0\}$. Note   that, by SA\ref{sass-homogeneity-f}-SA\ref{sass-homogeneity-ell},    $f(0,0)=0$ and $\ell(0,0)=0$, hence the optimal sequence at the origin is the 0 sequence of appropriate length, and $V_{d,\gamma_1,\gamma_2}(0)=0$ for any $d\in\Zp\cup\{\infty\}$, $\gamma_1\in\Rlp$ and $\gamma_2\in\R$. The extra condition imposed by Corollary \ref{cor-scaling-optimal-value-function} is on the homogeneity degree of $f$, which needs to be equal to $1$, or on the one of $\ell$ and $\jmath$, which needs to be equal to $0$. The first case corresponds to  \cite[Proposition 2]{Grimm-et-al-tac2005}, up to the minor difference that  here the cost can have an infinite horizon. The second case, namely when the homogeneity degree of $\ell$ and $\jmath$ is $0$, is new  as far as we know, and is important as it covers minimum time control and maximum hands off control as particular cases, as explained in more detail in Section \ref{subsect:sass-ell-jmath}.

\begin{rem}
Theorem \ref{thm-scaling-optimal-value-function} may also be useful when the primary control objective is to constrain the state in a given set $\mathcal{S}\subset\R^{n_x}$. This requirement can be incorporated in the definition of the stage cost using the function $\delta_{\mathcal{S}}^{\infty}$, which is zero when the argument is in $\mathcal{S}$ and infinite elsewhere, as $\delta_{\mathcal{S}}$ satisfies the desired homogeneity property when $\mathcal{S}$ is homogeneous, as discussed in more details in Section \ref{subsect:sass-ell-jmath}. In this case, even though the sequence of optimal inputs obtained at $x$ can be scaled to minimize another cost function on $\mathcal{R}_r(x)$ according to Theorem \ref{thm-scaling-optimal-value-function}, the obtained sequence of inputs ensures the cost  is bounded and thus that the solution remain in $\mathcal{S}$ for all future times. Indeed, if the solution for the scaled state did not remain in $\mathcal{S}$, the scaled cost would be infinite from the cost incurred by $\delta_{\mathcal{S}}^{\infty}$ when the state is not in $\mathcal{S}$. Since the scaled cost is finite, the state is always in $\mathcal{S}$.  \mbox{}\hfill $\Box$
\end{rem}

% We discuss below two particular cases of Theorem \ref{thm-scaling-optimal-value-function}, which could be formalized as lemmas later. 
% When SA\ref{sass-homogeneity-f} holds with $\nu=1$, (\ref{eq-thm-scaling-optimal-value-function}) gives  
% \begin{equation}
% \begin{array}{rllll}
% V_{d,\varepsilon^{-\mu},1}(\lambda^{r}(\varepsilon)x) & = & \varepsilon^{-\mu}V_{d,1,1}(\lambda^{r}(\varepsilon)x).
% \end{array}
% \end{equation}
% Hence, in view of (\ref{eq-thm-scaling-optimal-value-function}),
% \begin{equation}
% \begin{array}{rllll}
% V_{d,1,1}(\lambda^{r}(\varepsilon)x) & = & e^{\mu} V_{d,1,1}(x),
% \end{array}
% \end{equation}
% which corresponds to the result stated in \cite[Proposition 2]{Grimm-et-al-tac2005} with no terminal cost.

% When SA\ref{sass-homogeneity-ell} holds with $\mu=0$, (\ref{eq-thm-scaling-optimal-value-function}) gives  
% \begin{equation}
% \begin{array}{rllll}
% V_{d,1,\nu}(\lambda^{r}(\varepsilon)x) & = & V_{d,1,1}(\lambda^{r}(\varepsilon)x).
% \end{array}
% \end{equation}
% Hence, in view of (\ref{eq-thm-scaling-optimal-value-function}),
% \begin{equation}
% \begin{array}{rllll}
% V_{d,1,1}(\lambda^{r}(\varepsilon)x) & = & V_{d,1,1}(x).
% \end{array}
% \end{equation}
% \romain{This result may be useful for minimum time optimal control. Indeed, in this  case, when the target is the origin, we have $\ell(x,u)=\bm{1}_{0}(x)$ where $\bm{1}_{0}$ is the indicator function of the origin. Then, $\ell$ is homogeneous of degree $\mu=0$ with respect to any dilation pair according to SA\ref{sass-homogeneity-ell}. }

%%%%%%%%%%%%%%%%%%%%
\vspace{-0.1em}\section{Approximate VI scheme}\label{sect:dp}
%%%%%%%%%%%%%%%%%%%

We explain in this section how the scaling property of the optimal value function established in Theorem \ref{thm-scaling-optimal-value-function} can be exploited to derive a new approximation scheme for VI \cite{Bertsekas-book12(adp)} for homogeneous systems and homogeneous   costs satisfying SA\ref{sass-homogeneity-f}--SA\ref{sass-existence-optimal-value-function}.

%%%%%%%%%%%%
\vspace{-0.2em}\subsection{Value iteration} \label{subsect:vi}
%%%%%%%%%%%%%
We recall that VI is an iterative procedure to solve optimal control problems \cite{Bertsekas-book12(adp)}. We consider in this section an infinite-horizon cost, i.e., $d=\infty$ and $\jmath=0$ in {\eqref{eq-cost}}, and the goal is to derive (near-)optimal sequence of inputs for any state $x\in\R^{n_x}$. 
VI works as follows. Given an initial value function denoted $V_{0}:\R^{n_x}\to\Rlo$, at iteration $i+1\in\Zo$, an estimate of the optimal value function $V_{\infty,1,1}$ is obtained by solving, for any $x\in\R^{n_x}$,
\begin{equation}
V_{i+1}(x) = \dst\min_{u}\left[\ell(x,u)+V_{i}(f(x,u))\right],
\label{eq:vi}
\end{equation}
assuming the minimum exists. When solving (\ref{eq:vi}), we can retrieve a minimizer $h_{i+1}(x)$ for any $x\in\R^{n_x}$, which is used for control. Note that VI trivially accepts the finite-horizon case by initializing with $V_{0}=\jmath$ and stopping at $i=d$, but we do not consider this case here for convenience.

To exactly solve (\ref{eq:vi}) for any $x\in\R^{n_x}$ is in general impossible for system (\ref{eq-plant}) and stage cost $\ell$. Instead, most of the time we compute an approximation of $V_{i}$ on a compact set of interest, see, e.g., \cite{Bertsekas-18,busoniu-et-al-aut10}. The drawbacks of this approach are that (i) the selection of the compact set is not always rigorously established, (ii) the approximation is only defined on the considered set, (iii) we therefore need to make sure this set is forward invariant for the considered controlled system with the obtained control policy, which is not true in general, (iv) the computational cost may still be prohibitive. In the next section, we propose to exploit homogeneity to mitigate these issues.
%%%%%%%
\vspace{-0.2em}\subsection{Algorithm}\label{subsect:homVI}
%%%%%%%
The idea is to solve (\ref{eq:vi}) only on a given compact set $\mathcal{D}$ and to derive the value of $V_{i}$ elsewhere in the state space thanks to Theorem \ref{thm-scaling-optimal-value-function}, under SA\ref{sass-homogeneity-f}-SA\ref{sass-existence-optimal-value-function}. The proposed algorithm works as follows. 

The compact set $\mathcal{D}\subset\R^{n_x}$ is selected such that
\begin{equation}
\dst\bigcup_{x\in\mathcal{D}}\mathcal{R}_r(x)  =  \R^{n_x}\setminus\{0\},
\label{eq:mathcal-D-property}
\end{equation}
that is, the union of all  possible dilations of $x\in\mathcal{D}$ spans $\R^{n_x}\setminus\{0\}$. Examples of set $\mathcal{D}$ include the unit $n_x$-sphere or $\left\{x\in\R^{n_x}\,:\,|x_1|^{\frac{1}{\mu r_1}}+\ldots+|x_{n_x}|^{\frac{1}{\mu r_{n_x}}}=1\right\}$, where $(r_1,\allowbreak \ldots, r_{n_x})$ come from SA\ref{sass-homogeneity-f} and $\mu>0$ can take any value. Note that both of these examples are manifolds of degree $n_x-1$ while the state space is of dimension $n_x$. Again, the origin is excluded in (\ref{eq:mathcal-D-property}) as we know that the optimal value function is equal to $0$ in this case, as so are $f(0,0)$ and $\ell(0,0)$ in view of  SA\ref{sass-homogeneity-f}-SA\ref{sass-homogeneity-ell}.

We initialize the modified version of VI with $\underline{V}_{0}=\overline{V}_{0}=V_0$ where $V_0:\R^{n_x}\to\Rlo$ is homogeneous of degree $\mu$ with respect to dilation pair $(\lambda^{r},\lambda^{q})$, i.e., for any $x\in\R^{n_x}$ and any $\varepsilon>0$, $V_0(\lambda^{r}(\varepsilon)x)=\varepsilon^{\mu}V_0(x)$. We can simply take $V_0=0$ for instance. At iteration $i+1\in\Zo$, we generate a lower and an upper bound of $V_{i+1}$ in (\ref{eq:vi}), which are respectively denoted $\underline{V}_{i+1}$ and $\overline V_{i+1}$, like in relaxed dynamic programming \cite{Lincoln-Rantzer-tac06}. In particular, $\underline{V}_i(0)=\overline{V}_i(0)=0$ and, for any $x\in\mathcal{D}$ and $\varepsilon > 0$, 
\begin{equation}
\left\{\begin{aligned} 
\underline{V}_{i+1}(x)  &{}=   \dst\min_{u}\left[\ell(x,u)+\underline{V}_{i}(f(x,u))\right]\\
\underline{V}_{i+1}(\lambda^{r}(\varepsilon)x)   &{}=  \min\{\varepsilon^{\mu},\varepsilon^{\mu\nu^{i+1}}\} \underline{V}_{i+1}(x)
\end{aligned}\right.\label{eq:underline-vi}
\end{equation}
and 
\begin{equation}
\left\{\begin{aligned}
\overline{V}_{i+1}(x) &{}= \dst\min_{u}\left[\ell(x,u)+\overline{V}_{i}(f(x,u))\right]\\
\overline{V}_{i+1}(\lambda^{r}(\varepsilon)x) &{}= \max\{\varepsilon^{\mu},\varepsilon^{\mu\nu^{i+1}}\} \overline{V}_{i+1}(x).
\end{aligned}\right.\label{eq:overline-vi}
\end{equation}
In both (\ref{eq:underline-vi}) and (\ref{eq:overline-vi}), VI equation (\ref{eq:vi}) is only solved for $x$ in $\mathcal{D}$. The values of $\underline{V}_i$ and $\overline{V}_i$ when $x\notin\mathcal{D}$ do not require solving (\ref{eq:vi}): these are instead derived through the second equation of (\ref{eq:underline-vi}) and (\ref{eq:overline-vi}), respectively. Note that the two equations in (\ref{eq:underline-vi}) and (\ref{eq:overline-vi}) match when $\varepsilon=1$ as $\lambda^{r}(1)x=x$. 

\begin{rem} Compared to classical VI \eqref{eq:vi}, we only have to solve the first line of  \eqref{eq:underline-vi} and \eqref{eq:overline-vi} in the given set $\mathcal{D}$, e.g., in an  $n_x$-sphere {centered at the origin}. Moreover, an  implementation of \eqref{eq:vi} in a \emph{forward invariant} compact set will have the value functions  defined only for states in such set, which in contrast  the approximation scheme  \eqref{eq:underline-vi} and \eqref{eq:overline-vi} provide    values for the whole state space $\R^{n_x}$, in view of  their respective  second line.  However, executing the minimization   exactly for all states   in set $\mathcal{D}$ is still difficult in general,   we will study the effect of the induced approximation errors in future work.

\mbox{}\hfill$\Box$\end{rem} 
%\begin{rem} \diff{Again, the interest of \eqref{eq:} is to solve \eqref{eq:underline-vi} and \eqref{eq:overline-vi} in the given set $\mathcal{D}$, e.g., an $n_x$-sphere, and estimate the value function elsewhere.  Compared to classical VI, this set $\mathcal{D}$ is well defined and we can derive the value function in the whole state space.  However, executing VI exactly for this given set  is still difficult in general and we study approximation errors in future work.  \mbox{}\hfill }$\Box$\end{rem}  

%Compare to clarif we do not fix a compact set and do not derive its value only for this compact, here has to be selected as in 23 and we only have to solve VI in this set to derive its value on whole state space. Still, its hard to solve exactly the first equations in 24 25 and approximation scheme are needed as we do in section bla.

%%%%%%%%
\vspace{-0.2em}\subsection{Guarantees}\label{subsect:homVIanalysis}
%%%%%%%%
The next proposition\iflong\mbox{ }\else, whose proof is available in \cite{missing},\mbox{ }\fi ensures that $\underline{V}_i$ and $\overline{V}_i$ are indeed lower and upper bounds of $V_i$ on $\R^{n_x}$.
\begin{prop}\label{prop-vi} For any $i\in\Zo$ and any $x\in\R^{n_x}$,
\begin{equation}
\thickmuskip=15mu
\underline{V}_{i}(x)  \leq  V_i(x)  \leq  \overline{V}_i(x),
\label{eq:prop-vi}
\end{equation}
where $V_i$ is defined by iteration (\ref{eq:vi}) and is initialized with $\underline{V}_0(x)=V_0(x)=\overline{V}_0(x)$. In addition, when $\mu=0$ or $\nu=1$, $\underline{V}_i(x)=V_i(x)=\overline{V}_i(x)$ for any $x\in\R^{n_x}$. \hfill $\Box$
\end{prop}
%%%%%%%%%%%%%%
\iflong
\noindent\textbf{Proof:} In view of (\ref{eq:underline-vi}) and (\ref{eq:overline-vi}), $\underline{V}_i(0)=\overline{V}_i(0)=0$. Since SA\ref{sass-homogeneity-f}-SA\ref{sass-homogeneity-ell} hold, $V_i(0)=0$. Hence, we only need to prove that (\ref{eq:prop-vi}) holds for any $x\in\R^{n_x}\setminus\{0\}$, which is equivalent to proving that, for any $x\in\mathcal{D}$ and $\varepsilon>0$, 
\begin{equation}
\thickmuskip=15mu
\underline{V}_{i}(\lambda^{r}(\varepsilon)x)  \leq  V_i(\lambda^{r}(\varepsilon)x)  \leq  \overline{V}_i(\lambda^{r}(\varepsilon)x),
\label{eq:prop-vi-proof}
\end{equation}
in view of the property of $\mathcal{D}$ in (\ref{eq:mathcal-D-property}). We proceed by induction to prove (\ref{eq:prop-vi-proof}). 

\emph{$i=0$}. Property (\ref{eq:prop-vi-proof}) holds as $\underline{V}_0=\overline{V}_0=V_0$. 

\emph{$i\Rightarrow i+1$}. We assume that (\ref{eq:prop-vi-proof}) holds for $i\in\Zo$ and we aim at proving that this   is also satisfied for $i+1$. Let $x\in\mathcal{D}$ and $\varepsilon>0$. 
When $\varepsilon=1$, in view of (\ref{eq:underline-vi}),
\begin{equation}
\underline{V}_{i+1}(x)  =  \dst\min_{u}\left[\ell(x,u)+\underline{V}_i(f(x,u))\right].
\end{equation}
By assumption $\underline{V}_i \leq V_i$, thus, in view of (\ref{eq:vi}),
\begin{equation}
\begin{split}
\underline{V}_{i+1}(x) & {}\leq  \dst\min_{u}\left[\ell(x,u)+V_i(f(x,u))\right]\\
& {}=  V_{i+1}(x).
\end{split}
\end{equation}
When $\varepsilon\neq 1$, from (\ref{eq:vi}),
\begin{align}
V_{i+1}(\lambda^{r}(\varepsilon)x) & {}= \dst\min_{u}\left[\ell(\lambda^{r}(\varepsilon)x,u)+V_i(f(\lambda^{r}(\varepsilon)x,u))\right]\nonumber\\
& {}= \dst\min_{u}\left[\ell(\lambda^{r}(\varepsilon)x,\lambda^{q}(\varepsilon)\lambda^{q}(\varepsilon^{-1})u)\right.\\&\qquad\left.{}+V_i(f(\lambda^{r}(\varepsilon)x,\lambda^{q}(\varepsilon)\lambda^{q}(\varepsilon^{-1})u))\right].\nonumber
\end{align}
We derive from SA\ref{sass-homogeneity-f}-\ref{sass-homogeneity-ell} and the fact that $\lambda^{q}(\varepsilon^{-1})\R^{n_u}=\R^{n_u}$,
\begin{align}
\MoveEqLeft V_{i+1}(\lambda^{r}(\varepsilon)x) \nonumber \\  &{}= \dst\min_{u}\left[\varepsilon^{\mu}\ell(x,\lambda^{q}(\varepsilon^{-1})u) + V_i(\lambda^{r}(\varepsilon)^{\nu}f(x,\lambda^{q}(\varepsilon^{-1})u))\right]\nonumber\\
&{}= \dst\min_{v}\left[\varepsilon^{\mu}\ell(x,v)+V_i(\lambda^{r}(\varepsilon)^{\nu}f(x,v))\right].
\label{eq:proof-prop-V-i+1}
\end{align}
Since $\underline{V}_i\leq V_i$ by assumption,
\begin{equation}
V_{i+1}(\lambda^{r}(\varepsilon)x)  \geq  \dst\min_{v}\left[\varepsilon^{\mu}\ell(x,v)+\underline{V}_i(\lambda^{r}(\varepsilon)^{\nu}f(x,v))\right].
\end{equation}
We have from (\ref{eq:underline-vi}) that $\underline{V}_i(\lambda^{r}(\varepsilon)^{\nu}f(x,v))=\min\{\varepsilon^{\mu\nu},\varepsilon^{\mu\nu^{i+1}}\}\underline{V}_i(f(x,v))$ for any $v\in\R^{n_u}$, hence
\begin{align}
\MoveEqLeft V_{i+1}(\lambda^{r}(\varepsilon)x) \nonumber\\& {}\geq  \dst\min_{v}\left[\varepsilon^{\mu}\ell(x,v)+\min\{\varepsilon^{\mu\nu},\varepsilon^{\mu\nu^{i+1}}\}\underline{V}_i(f(x,v))\right]\nonumber\\
& {}\geq  \dst\min\{\varepsilon^{\mu},\varepsilon^{\mu\nu},\varepsilon^{\mu\nu^{i+1}}\}\min_{v}\left[\ell(x,v)+\underline{V}_i(f(x,v))\right].
\end{align}
Since $\min\{\varepsilon^{\mu},\varepsilon^{\mu\nu},\varepsilon^{\mu\nu^{i+1}}\}=\min\{\varepsilon^{\mu},\varepsilon^{\mu\nu^{i+1}}\}$ and in view of (\ref{eq:underline-vi}),
\begin{equation}
\begin{split}
V_{i+1}(\lambda^{r}(\varepsilon)x) &  {} \geq  \dst\min\{\varepsilon^{\mu},\varepsilon^{\mu\nu^{i+1}}\}\underline{V}_{i+1}(x)\\
& {} =  \underline{V}_{i+1}(\lambda^{r}(\varepsilon)x).
\end{split}
\end{equation}
The proof that  $V_{i+1}(\lambda^{r}(\varepsilon)x)\leq\overline{V}_{i+1}(\lambda^{r}(\varepsilon)x)$ follows similar lines as above. Consequently, (\ref{eq:prop-vi-proof}) holds for $i+1$. The proof by induction is complete: (\ref{eq:prop-vi}) is guaranteed. 

When $\mu=0$ or $\nu=1$, $\varepsilon^{\mu}=\varepsilon^{\mu\nu^{i}}$ for any $i\in\Zo$. Hence, since $\underline{V}_0=\overline{V}_0=V_0$, we derive from  (\ref{eq:underline-vi}) and (\ref{eq:overline-vi}) that $\underline{V}_i(x)=\overline{V}_i(x)$ for any $i\in\Zo$ and $x\in\R^{n_x}$. Consequently, $\underline{V}_i(x)=\overline{V}_i(x)=V_i(x)$ for any $i\in\Zo$ and $x\in\R^{n_x}$ in view of (\ref{eq:prop-vi}). \hfill $\blacksquare$
\fi
%%%%%%%%%%%%%%%%
% 

Proposition \ref{prop-vi} ensures that the exact value function  given by VI at any iteration $i\in\Zo$, i.e., $V_i$ in (\ref{eq:vi}) is lower and upper-bounded by $\underline{V}_i$ and $\overline{V}_i$ as defined in (\ref{eq:underline-vi}) and (\ref{eq:overline-vi}), respectively. We can then exploit these bounds to construct a control law by taking, at any iteration $i\in\Zo$ and for any $x\in\R^{n_x}$,
\begin{equation}
\thickmuskip=15mu
\hat{h}_{i}(x) \in  \dst\min_{u}\left[\ell(x,u)+\widehat{V}_{i}(f(x,u))\right],
\label{eq:hat-h-i}    
\end{equation}
where $\widehat{V}_i:\R^{n_x}\to\R$ is any function such that $\underline{V}_i \leq \widehat{V}_i \leq \overline{V}_i$. Compared to relaxed dynamic programming  \cite{Lincoln-Rantzer-tac06} where lower and upper-bounds on $V_i$ are used to derive near-optimal inputs, $\underline{V}_i$ and $\overline{V}_i$ in (\ref{eq:underline-vi}) and (\ref{eq:overline-vi}) do not require to solve minimization problems on the whole state space, but only on $\mathcal{D}$, which is favourable for practical implementation. Moreover, when $\mu=0$ or $\nu=1$ as in Corollary \ref{cor-scaling-optimal-value-function}, we  only need to solve (\ref{eq:vi}) on $\mathcal{D}$ at each iteration to deduce the value of $V_i$ on the whole state space $\R^{n_x}$, in view of the last part of Proposition \ref{prop-vi}. We investigate numerically the tightness of the bounds of the proposed scheme in Section \ref{sect:ex}.

%\romain{It would be interesting to see whether we could prove whether $\underline{V}_i$ and $\overline{V}_i$ are closed in a neighborhood of $\mathcal{D}$, but let see first what the simulations give.}

%\romain{It is possible to derive a tighter upper-bound $\overline{V}_i$ on $V_i$, which is more difficult to implement. We will see whether we need to mention it.} 

%{An obvious candidate for $\widehat V_i$ is $\frac{\overline V_i + \underline V_i}{2}$. The gap between $V_i$ and $\widehat V_i$ is also $\overline V_i - \underline V_i$. A priori, the errors for $\lambda^r(\varepsilon)x$ should be enough to make $\overline V_i$ diverge from $\underline V_i$ \emph{everywhere}, in particular for $i>>0$. I am curious to see if we can use $\widehat V_i$ instead of  $\overline V_i$  and  $\underline V_i$ in the minimization steps of \eqref{eq:underline-vi} and \eqref{eq:overline-vi}, as an approximation to the true $V_i$ to minimize divergence.}
% Bounds match when $\varepsilon=1$ ? 

% A natural question is whether $\underline{V}_i$ and $\overline{V}_i$ are tight bounds of $V_i$. The next result states that this is the case when the state is closed to $\mathcal{D}$, which is not surprizing in view of (\ref{eq:underline-vi}) and (\ref{eq:overline-vi}).

% \begin{prop} Result on what happens when $\varepsilon\to 1$ XXX.
% \end{prop}

%%%%%%%%
\vspace{-0.1em}\section{Ensuring the standing assumptions}\label{sect-SA} %Sufficient conditions that guarantee the standing assumptions
%%%%%%%%
In this section, we provide examples of vector fields $f$ in (\ref{eq-plant}) and stage cost $\ell$, which verify SA\ref{sass-homogeneity-f}-\ref{sass-homogeneity-ell}.   SA\ref{sass-existence-optimal-value-function}, on the other hand, can be analysed using the results in \cite{Keerthi-Gilbert-tac85} as already mentioned.\iflong \else\mbox{ }The proofs of this section can be found in \cite{missing}.\fi 

\vspace{-0.2em}\subsection{Standing Assumption \ref{sass-homogeneity-f}}\label{subsect-sass-f}

The next result may be used to ensure the satisfaction of SA\ref{sass-homogeneity-f}.

\begin{prop}\label{prop-sass-f} Suppose that there exist $r=(r_1,\ldots,r_{n_x})\in\R^{n_x}_{>0}$, $q=(q_1,\ldots,q_{n_u})\in\R^{n_u}_{>0}$, scalar fields $g_{i,j}:\R^{n_x}\times\R^{n_u}\to\R$ that are homogeneous of degree $\nu_{i,j}$ with respect to a dilation pair $(\lambda^r,\lambda^q)$ and $N\in\Zp$ such that, for any $i\in\{1,\ldots,n_x\}$ and $(x,u)\in\R^{n_x}\times\R^{n_u}$,
\begin{equation}
f_i(x,u)  =  \dst\sum_{j=0}^{N} g_{i,j}(x,u).
\label{eq-prop-fi}
\end{equation}
For any $(x,u,w)\in\R^{n_x}\times\R^{n_u}\times\R$, the vector field $\tilde{f}=(\tilde{f}_1,\ldots,\tilde{f}_{n_x},\tilde{f}_{n_x+1})$  with
\begin{equation}
\begin{split}
\tilde{f}_i(x,w,u)   & {} =   \dst\sum_{j=0}^{N} w^{(r_i \nu -\nu_{i,j})\nu}g_{i,j}(x,u)\\
\tilde{f}_{n_x+1}(x,w,u)  & {} = \sabs{w}^{\nu},
\end{split}
\end{equation}
and $\nu:=\dst 1+\max\left\{\frac{\nu_{i,j}}{r_i}\,:\,(i,j) \in\{1,\ldots,n_x\}\times\{1,\ldots,N\}\right\}\allowbreak\in\Zp$, 
satisfies SA\ref{sass-homogeneity-f}, in particular it is homogeneous of degree $\nu\in\Z$ with respect to the dilation pair $(\lambda^{(r,1/\nu)},\lambda^q)$.

\mbox{}\hfill $\Box$ \end{prop}

\iflong
\noindent\textbf{Proof:} Let $(x,u,w)\in\R^{n_x}\times\R^{n_u}\times\R$ and $\varepsilon>0$. From the homogeneity property of $g_{i,j}$ for any $(i,j)\in\{1,\ldots,n_x\}\times\{1,\ldots,N\}$,
\begin{equation}
\begin{split}
\MoveEqLeft \tilde{f}_i(\lambda^{r}(\varepsilon)x,\lambda^{1/\nu}(\varepsilon)w,\lambda^{q}(\varepsilon)u) \\
&{}= \dst\sum_{j=0}^{N} \varepsilon^{r_i \nu -\nu_{i,j}} w^{(r_i \nu -\nu_{i,j})\nu}\varepsilon^{\nu_{i,j}}g_{i,j}(x,u) \\
&{}= \varepsilon^{r_i \nu} \dst\sum_{j=0}^{N} w^{(r_i \nu -\nu_{i,j})\nu}g_{i,j}(x,u)\\
&{}= \varepsilon^{r_i \nu} \tilde{f}_i(x,w,u).
\end{split}
\end{equation}
On the other hand, $\tilde{f}_{n_x+1}(\lambda^{r}(\varepsilon)x,\lambda^{1/\nu}(\varepsilon)w,\lambda^{q}(\varepsilon)u)=\sabs{\varepsilon^{1/\nu}w}^\nu = \varepsilon \sabs{w}^{\nu}=\varepsilon^{\nu\frac{1}{\nu}}\tilde{f}_{n_x+1}(x,w,u)$. Hence, $\tilde f$ is homogeneous of degree $\nu\in\Z$ with respect to dilation pair $(\lambda^{(r,1/\nu)},\lambda^q)$, which concludes the proof. \hfill $\blacksquare$
\fi

Proposition \ref{prop-sass-f} implies that when the components of $f$, namely the $f_i$'s, can be written as the sum of (scalar) homogeneous functions with respect to the same dilation pair, SA\ref{sass-homogeneity-f} can always be ensured by adding a variable $w$, whose dynamics is $w^{+}=\sabs{w}^{\nu}$. By initializing $w$ at $1$, we note that $w(k)=1$ for any $k\in\Zo$ and we have  $f(x,u)=(\tilde{f}_1(x,1,u),\ldots,\tilde{f}_{n_x}(x,1,u))$. The idea to add the extra dynamics $w$ to ensure the desired homogeneity property is common in the homogeneous systems literature, see for instance \cite{baillieul-na80} where a different homogeneity property is studied for polynomial vector fields.  Notice that there is no loss of generality  by considering the same integer $N$ for all components $\tilde{f}_i$ of $f$ in (\ref{eq-prop-fi}) as we can always add zero-valued functions if needed.

Proposition \ref{prop-sass-f} applies to polynomial vector fields as a particular case. Indeed, suppose $f$ is a polynomial vector field of degree $d\in\Zo$, then $f_i$ can be written as in (\ref{eq-prop-fi}) with $N= \left(\begin{array}{c}n_x+n_u+d \\ d\end{array}\right)-1$ and $g_{i,j}$ are monomials of $(x,u)$ multiplied by a constant. Proposition \ref{prop-sass-f} also covers the case where $f$ is a rational function of components of $(x,u)$ using the same technique as in \cite{baillieul-na80}.   

Proposition \ref{prop-sass-f} reduces the homogeneity analysis of $f$ to those of the components of $f_i$, which may be derived from the next simple though useful lemma. Its proof directly follows from the homogeneity properties of the considered functions and is therefore omitted.

\begin{lem}\label{lem-homogeneity} Consider homogeneous scalar functions $g_1,g_2:\R\to\R$ of degree $\nu_1$ and $\nu_2$ respectively, then
\begin{enumerate}
\item[(i)] for any $\lambda\in\R$, $\lambda g_1$ is homogeneous of degree $\nu_1$,
\item[(ii)] $g_1\circ g_2$ is homogeneous of degree~$\nu_1+\nu_2$,
\item[(iii)] $(z_1,z_2)\mapsto g_1(z_1)g_2(z_2)$ is homogeneous of degree~$\nu_1+\nu_2$,
\item[(iv)] $(z_1,z_2,w)\mapsto w^{\max(\nu_1,\nu_2)}\left(w^{-\nu_1}g_1(z_1)+w^{-\nu_2}g_2(z_2)\right)$ is homogeneous of degree $\max(\nu_1,\nu_2)$. \hfill $\Box$
\end{enumerate}
\end{lem}

Examples of  functions $g_1$ and $g_2$ as in Lemma \ref{lem-homogeneity} include $z\mapsto z^{a}$ with $a\in\Zp$, $z\mapsto \sabs{z}^{a}$ with $a\in\Rlp$, ReLU functions $z\mapsto \max(0,z)$ and $z\mapsto \max(z,\epsilon z)$ with $\epsilon\in(0,1)$ to mention a few. In view of Lemma \ref{lem-homogeneity}, any vector field $f=(f_1,\ldots,f_{n_x})$ which involves compositions, products or linear combinations of such elementary homogeneous functions can be suitably modified using variable $w$ as in Proposition \ref{prop-sass-f} to ensure the satisfaction of SA\ref{sass-homogeneity-f}.

\vspace{-0.2em}\subsection{Standing Assumption \ref{sass-homogeneity-ell}}\label{subsect:sass-ell-jmath}

First, we have the next lemma, which shows that the linear combination of homogeneous stage (terminal) costs of different degrees but with respect to the same dilation pair can be modified by adding an extra dummy variable $w$ as in Proposition \ref{prop-sass-f} to ensure the satisfaction of SA\ref{sass-homogeneity-ell}.

\begin{lem}\label{lem-sass-ell-linear-combination} Let $m\in\Zp$ and $\ell_1,\ldots,\ell_m:\R^{n_x}\times\R^{n_u}\to\overline{\R}_{\geq 0}$ (respectively, $\jmath_1,\ldots,\jmath_m:\R^{n_x}\to\overline{\R}_{\geq 0}$) be  homogeneous with respect to a dilation pair $(\lambda^{r},\lambda^{q})$ (respectively, $\lambda^r$) of degree $\mu_{1},\ldots,\mu_m\in\Zo$. Let $\mu=\max\{\mu_i\,:\,i\in\{1,\ldots,m\}\}$. Then, for any $(x,w,u)\in\R^{n_x+1}\times\R^{n_u}$, 
$\ell(x,w,u) = \dst\sum_{i=1}^{m}w^{\mu-\mu_i}\ell_i(x,u)$  
is homogeneous of degree $\mu$ with respect to the dilation pair $(\lambda^{(r,1)},\lambda^{q})$ (respectively, $\jmath(x,w,u)=\dst \dst\sum_{i=1}^{m}w^{\mu}\jmath_i(x)$ is homogeneous of degree $\mu$ with respect to $\lambda^{(r,1)}$). When $\mu_1=\ldots=\mu_m=\mu$, $\ell=\dst\sum_{i=1}^{m}\ell_i$ is homogeneous of degree $\mu$ with respect to the dilation  pair $(\lambda^{r},\lambda^{q})$ (respectively, $\jmath=\dst \dst\sum_{i=1}^{m}\jmath_i$ is homogeneous of degree $\mu$ with respect to $\lambda^r$). \hfill $\Box$
\end{lem}

\iflong
\noindent\textbf{Proof:} Let $(x,w,u)\in\R^{n_x+1}\times\R^{n_u}$ and $\varepsilon>0$. In view of the homogeneity properties of $\ell_1,\ldots,\ell_m$,
\begin{align}
\ell(\lambda^{r}(\varepsilon)x,\lambda^{1}(\varepsilon)w,\lambda^{q}(\varepsilon)u) & {} = \dst\sum_{i=1}^{m}\varepsilon^{\mu-\mu_i}w^{\mu-\mu_i}\varepsilon^{\mu_i}\ell_i(x,u) \nonumber\\
& {} = \dst\varepsilon^{\mu}\ell(x,w,u). 
\end{align}
We have proved that $\ell$ is homogeneous of degree $\mu$ with respect to dilation pair $(\lambda^{(r,1)},\lambda^{q})$. When $\mu_1=\ldots=\mu_m=\mu$, we directly derive the desired result from the above development. Finally, the case of $\jmath$ similarly follow. \hfill $\blacksquare$
\fi

%The next lemma provides examples of quadratic(-like) stage costs verifying SA\ref{sass-homogeneity-ell}.
% dillation map is given from SA1 and we want to know in which SA2 is verified with quadratic costs.
The next lemma provides a condition for the dilation maps  from SA\ref{sass-homogeneity-f} under which  SA\ref{sass-homogeneity-ell} is verified with quadratic(-like) costs.

\begin{lem}\label{lem-sass-ell} The following holds.
\begin{enumerate}
\item[(i)] Let $\ell:(x,u) \mapsto |x|_Q + |u|_R$, where $Q=Q^{\top}\geq 0$, 
$R=R^{\top}\geq 0$. For any $\mu\in\Zp$, SA\ref{sass-homogeneity-ell} holds with degree $2\mu$ and for any dilation pair $(\lambda^{r},\lambda^{q})$ of the form  $r=\mu(1,\ldots, 1)\in\R^{n_x}$ and $q=\mu(1,\ldots,1)\in\R^{n_u}$.
% \item[(ii)] Let $\ell:(x,u) \mapsto x^\top Qx + u^\top R u$, where $Q=\text{diag}(\theta_1,\ldots,\theta_{n_x})$, 
% $R=\text{diag}(\vartheta_1,\ldots,\vartheta_{n_u})$, $\theta_i,\vartheta_j\geq 0$ for $i\in\{1,\ldots,n_x\}$ and $j\in\{1,\ldots,n_u\}$. For any $\mu\in\Zp$, SA\ref{sass-homogeneity-ell} holds with degree $2\mu$ and for  any dilation pair $(\lambda^{r},\lambda^{q})$ of the form  $r=(r_1,\ldots,r_{n_x})$ and $q=(q_1,\ldots,q_{n_u})$ with  $r_i=p$ and $q_j=p$ for any $(i,j)\in\{1,\ldots,n_x\}\times\{1,\ldots,n_u\}$ such that $\theta_i\neq 0$ and $\vartheta_j\neq 0$.
\item[(ii)] For any $r=(r_1,\ldots,r_{n_x})\in\R^{n_x}_{>0}$, $q=(q_1,\ldots,\allowbreak q_{n_u})\in\R^{n_u}_{>0}$, $\mu\in\Zp$, $Q=Q^{\top}\geq 0$,  and $R=R^{\top}\geq 0$, $\dst\ell:(x,u)\mapsto \left|\left(\sabs{x_1}^{\frac{\mu}{2r_1}},\ldots,\sabs{x_{n_x}}^{\frac{\mu}{2r_{n_x}}}\right)\right|_Q +  \left| \left(\sabs{u_1}^{\frac{\mu}{2q_1}},\ldots,\sabs{u_{n_u}}^{\frac{\mu}{2r_{n_u}}}\right)\right|_R$ satisfies SA\ref{sass-homogeneity-ell} with degree $\mu$ and  dilation pair $(\lambda^{r},\lambda^{q})$. \hfill $\Box$
\end{enumerate}
\end{lem}

\iflong
\noindent\textbf{Proof:} Let $(x,u)\in\R^{n_x}\times\R^{n_u}$ and $\varepsilon>0$. 

\emph{Item (i):} From the definitions of $\ell$, $r$ and $q$, we have
\begin{align}
\ell(\lambda^{r}(\varepsilon)x,\lambda^{q}(\varepsilon)u) & {} =  \left|(\varepsilon^\mu x)\right|_Q + \left|(\varepsilon^\mu u)\right|_R\nonumber\\
& {} =  \dst \varepsilon^{2\mu}\ell(x,u),
\end{align}
from which we derive that $\ell$ is homogeneous of degree $2\mu$ with respect to dilation pair $(\lambda^{p},\lambda^{q})$.

\emph{Item (ii):} From the definitions of $\ell$, $r$ and $q$, we have
 \begin{equation}
 \begin{array}{rlll}
 \ell(\lambda^{r}(\varepsilon)x,\lambda^{q}(\varepsilon)u) & = & \dst\sum_{i=1}^{n_x} \theta_i\varepsilon^{2\mu}x_i^{2} + \sum_{j=1}^{n_u} \vartheta_j\varepsilon^{2\mu}u_j^{2}\\
 & = & \dst \varepsilon^{2\mu}\ell(x,u),
 \end{array}
 \end{equation}
 from which we derive that $\ell$ is homogeneous of degree $2\mu$ with respect to dilation pair $(\lambda^{p},\lambda^{q})$.

\emph{Item (ii):} Let $r=(r_1,\ldots,r_{n_x})\in\R^{n_x}_{>0}$, $q=(q_1,\ldots,q_{n_u})\in\R^{n_u}_{>0}$, $\mu\in\Zp$, $Q=Q^{\top}\geq 0$,  and $R=R^{\top}\geq 0$. Since, for any $z\in\R$ and $c_1,c_2\in\Rlp$, $\sabs{c_1 z_1 z_2}^{c_2}=c_1^{c_2}\sabs{z_1}^{c_2}$ in view of the definition of $\sabs{\cdot}$, we derive
\begin{equation}
\begin{split}
\MoveEqLeft \ell(\lambda^{r}(\varepsilon)x,\lambda^{q}(\varepsilon)u) \\ &{} =  \left|\left(\sabs{\varepsilon^{r_1}x_1}^{\frac{\mu}{2r_1}},\ldots,\sabs{\varepsilon^{r_{n_x}}x_{n_x}}^{\frac{\mu}{2r_{n_x}}}\right)\right|_Q \\
& \quad {} +  \left| \left(\sabs{\varepsilon^{q_1} u_1}^{\frac{\mu}{2q_1}},\ldots,\sabs{\varepsilon^{q_{n_u}} u_{n_u}}^{\frac{\mu}{2r_{n_u}}}\right)\right|_R \\
& {} =   \varepsilon^{\mu}\ell(x,u).
\end{split}
\end{equation}
Hence $\ell$ is homogeneous of degree $\mu$ with respect to dilation pair $(\lambda^r,\lambda^q)$. \hfill $\blacksquare$
\fi

Item (i) of Lemma \ref{lem-sass-ell} shows that quadratic stage costs satisfy SA\ref{sass-homogeneity-ell} when SA\ref{sass-homogeneity-f} holds with $r_1=\ldots=r_{n_x}=q_1=\ldots=q_{n_u}$. % Item (ii) of Lemma \ref{lem-sass-ell} presents a variation for the specific case where  the  weight matrices are diagonal \romain{(I feel we can remove item (ii))}. 
When the dilation pair does not satisfy this property, the stage cost can always be modified as in item (ii) of Lemma \ref{lem-sass-ell} to enforce SA\ref{sass-homogeneity-ell}. This change of cost was proposed in the paragraph before \cite[Proposition 1]{Grimm-et-al-tac2005}, except that the function $\sabs{\cdot}$ is used here instead of the absolute value for generality. Regarding the terminal cost $\jmath$, the results of Lemma \ref{lem-sass-ell} apply by setting $R=0$. 

We have seen at the end of Section \ref{sect-main-results} that having a stage cost and a terminal cost of  homogeneity degree $0$ can be exploited to simplify the computation of the optimal value functions and of a sequence of optimal inputs for any $x\in\R^{n_x}$. It appears that the examples of stage and terminal costs in Lemma \ref{lem-homogeneity} are not homogeneous of degree zero. The next lemma provides relevant examples ensuring homogeneity with degree $0$, but before that we need to recall the notion of   homogeneous sets.

\begin{defn} Let $r=(r_1,\ldots,r_{n_x})\in\R^{n_x}_{>0}$. A set $\mathcal{S}\subseteq\R^{n_x}$ is \emph{homogeneous with respect to a dilation $\lambda^{r}$} when, for any $x\in\R^{n_x}$, $x\in \mathcal{S}$ if and only is $\lambda^{r}(\varepsilon)x\in \mathcal{S}$.\hfill $\Box$
\end{defn}

We write for short that a set is homogeneous when there exists a dilation with respect to which it is homogeneous. Examples of homogeneous sets include  $\{0\}$ and $\text{vect}(e_{i_1},\ldots,e_{i_m})$ where $i_1,\ldots,i_m\in\{1,\ldots,n_x\}$ and  $e_1,\ldots,e_{n_x}$ are the standard unit vectors of $\R^{n_x}$ for any dilation map. We are ready to present examples of functions, which can be used to define stage and terminal costs ensuring SA\ref{sass-homogeneity-ell} with $0$-homogeneity degree. 

\begin{lem}\label{lem-sass-ell-degree-0} Let $t=(t_1,\ldots,t_{m})\in\R^{m}$ with $m\in\Zp$ and $\mathcal{S}\subseteq\R^{m}$ be a homogeneous set with respect to $\lambda^{t}$. For any $c\in\overline\R_{\geq 0}$,  $\delta^{c}_{\mathcal{S}}$ is homogeneous of degree $0$ with respect to dilation $\lambda^t$. \hfill $\Box$
\end{lem}

%%%%%%%%%%%%%%%%%%%%
\iflong
\noindent\textbf{Proof:} Let $y\in\R^{m}$, $c\in\overline\R_{\geq 0}$ and $\varepsilon>0$, $\delta^{c}_{\mathcal{S}}(\lambda^{t}(\varepsilon)y)=0$ when $\lambda^{t}(\varepsilon)y\in\mathcal{S}$, which is equivalent to $y\in\mathcal{S}$ by homogeneity of $\mathcal{S}$, and  $\delta^{c}_{\mathcal{S}}(\lambda^{t}(\varepsilon)y)=c$ when $y\not\in\mathcal{S}$ for the same reason. Hence $\delta^{c}_{\mathcal{S}}(\lambda^{t}(\varepsilon)y)=\delta^{c}_{\mathcal{S}}(y)$, which corresponds to the desired result. \hfill $\blacksquare$
\fi
%%%%%%%%%%%%%%%%%%

In view of Lemma \ref{lem-sass-ell-degree-0} and the fact that the linear combination of $0$-degree homogeneous functions is homogeneous of degree $0$ in view of Lemma \ref{lem-sass-ell-linear-combination}, we derive that important optimal control problems correspond to the stage and terminal costs of degree $0$ including 
\begin{itemize}
\item \emph{minimum time control} to a homogeneous set $\mathcal{S}_x\subset\R^{n_x}$, in which case  $d=\infty$, $\ell=\delta^{1}_{\mathcal{S}_x}$ and thus $\jmath=0$, %%%  
\item \emph{maximum hands off control} where the objective is to drive the state in a given homogeneous set $\mathcal{S}_x\subset\R^{n_x}$ in $d\in\Zp$ steps while using a zero input as frequently as possible, see \cite{nagahara-et-al-tac15,nagahara-et-al-max-hands-off-dt}. In this case, $\ell(x,u)=|u|_{0}$ for any $(x,u)\in\R^{n_x}\times\R^{n_u}$ and $\jmath=\delta^{\infty}_{\mathcal{S}_x}$; noting that the $\ell_0$ norm is also homogeneous of degree $0$ with respect to any dilation map.
\end{itemize}

% In view of Lemma \ref{lem-scaling-ell},
% \begin{equation}
% \begin{array}{rlllll}
% V_{d,\gamma}(\lambda^{r}(\varepsilon)x) & \leq & \dst\sum_{k=0}^{d}\gamma^{k}\varepsilon^{\mu\nu^k}\ell\left(\phi(k,x,\bm{u}|_{k}^{\star}(x)),u_k^{\star}(x)\right).
% \end{array}
% \end{equation}
% On the other hand, let $\bm{v}_{{d}}^{\star}(\lambda^{r}(\varepsilon)x)$ be a minimizing sequence of inputs for $J_{\gamma,d}(\lambda^{r}(\varepsilon)x,\bm{u})$, i.e. $V_{d,\gamma}(\lambda^{r}(\varepsilon)x)=J_{\gamma,d}(\lambda^{r}(\varepsilon)x,\bm{v}^{\star}(\lambda^{r}(\varepsilon)x))$. Hence
% \begin{equation}
% \begin{array}{rlllll}
% V_{d,\gamma}(\lambda^{r}(\varepsilon)x) & = & \dst\sum_{k=0}^{d}\gamma^{k}\ell\Big(\phi(k,\lambda^{r}(\varepsilon)x,\bm{v}^{\star}(\lambda^{r}(\varepsilon)x)|_k),v_k^{\star}(\lambda^{r}(\varepsilon)x)\Big) \\
% & = & \dst\sum_{k=0}^{d}\gamma^{k}\ell\Big(\phi(k,\lambda^{r}(\varepsilon)x,\Lambda^{q}_{{\nu,d}}(\varepsilon)\Lambda^{q}_{{\nu,d}}(\varepsilon^{-1})\bm{v}^{\star}(\lambda^{r}(\varepsilon)x)|_k),\lambda^{q}(\varepsilon)^{\nu^{k}}\lambda^{q}(\varepsilon^{-1})^{\nu^{k}}v_k^{\star}(\lambda^{r}(\varepsilon)x)\Big).
% \end{array}
% \end{equation}
% We deduce from Lemma \ref{lem-scaling-ell} that
% \begin{equation}
% \begin{array}{rlllll}
% V_{d,\gamma}(\lambda^{r}(\varepsilon)x) & = & \dst\sum_{k=0}^{d}\gamma^{k}\varepsilon^{\mu\nu^k}\ell\Big(\phi(k,x,\Lambda^{q}_{{\nu,d}}(\varepsilon^{-1})\bm{v}^{\star}(\lambda^{r}(\varepsilon)x)|_k),\lambda^{q}(\varepsilon^{-1})^{\nu^{k}}v_k^{\star}(\lambda^{r}(\varepsilon)x)\Big).
% \end{array}
% \end{equation}

%%%%%%%%%%%%%%%%%%%%
\vspace{-0.1em}\section{Numerical case study}\label{sect:ex}
We illustrate the results of Section \ref{sect-SA}. In particular, we   study how conservative the bounds in   
Proposition~\ref{prop-vi} are. 
\vspace{-0.2em}\subsection{System and cost}
We consider the optimal control of the discretized model of a van der Pol oscillator using Euler's direct method, $
(x_1^{+},x_2^{+}) =  ( x_1 + T x_2 , x_2 + T (a (1-x_1^2)x_2 - b x_1 + u))$
%\begin{equation}
%\begin{split}
%x_1^{+} & {} =  x_1 + T x_2 \\
%x_2^{+} & {} =  x_2 + T \left(a (1-x_1^2)x_2 - b x_1 + u \right),
%\label{eq:sys-ex-van-der-pol}
%\end{split}
%\end{equation}
where $x_1,x_2\in\R$ are the states, $u\in\R$ is the control input, $a,b\in\Rlp$ is a parameter and $T\in\Rlp$ is the sampling period used for discretization.

To ensure the satisfaction of SA\ref{sass-homogeneity-f}, we add a variable  $x_3$ so that the extended system is of the form of (\ref{eq-plant}) with $x:=(x_1,x_2,x_3)$ and $
  f(x,u) := 
(x_1 x_3^2 + T x_2 x_3^2, 
 x_2 x_3^2 +  T a (x_3^2-x_1^2)x_2- T b x_1 x_3^2 + T u,
x_3^{3})  
$.
%\begin{equation}
%f(x,u)  = 
%\left(\begin{array}{c}
%x_1 x_3^2 + T x_2 x_3^2 \\
%\left(\begin{aligned}x_2 x_3^2  &{}+  T a (x_3^2-x_1^2)x_2\\ 
%                        &{}-  T b x_1 x_3^2 + T u  \end{aligned}\right)\\
%x_3^{3}
%\end{array}\right).\label{eq:ex-f-van-der-Pol}
%\end{equation}
Like in Proposition~\ref{prop-sass-f}, by initializing $x_3(0)=1$, the $(x_1,x_2)$-component of the solutions to (\ref{eq-plant}) with $f$ defined above matches the solution to the original system initialized at the same point and,  with a careful choice of the stage cost in the sequel, will solve the same optimal control problem.
The vector field $f$  satisfies SA\ref{sass-homogeneity-f} with $r=(1,1,1)$, $q=3$ and $\nu=3$. Indeed, for any $x\in\R^{3}$, $u\in\R$ and  $\varepsilon>0$,
$
f(\lambda^{r}(\varepsilon)x,\lambda^{q}(\varepsilon)u) =  
(
\varepsilon^3 x_1 x_3^2 + T \varepsilon^3 x_2 x_3^2,
\varepsilon^{3}x_2 x_3^2  + T a\varepsilon^3 (x_3^2-x_1^2)x_2    - Tb \varepsilon^3 x_1 x_3^2+ \varepsilon^{3} Tu   ,
\varepsilon^{3} x_3^{3})  = \lambda^{r}(\varepsilon)^3 f(x,u).
$
We  fix $a=b=T=1$ in the sequel.
%\begin{align}
%\MoveEqLeft f(\lambda^{r}(\varepsilon)x,\lambda^{q}(\varepsilon)u)\nonumber\\&{}=  
%\left(\begin{array}{c}
%\varepsilon^3 x_1 x_3^2 + T \varepsilon^3 x_2 x_3^2 \\
%\left(\begin{aligned}\varepsilon^{3}x_2 x_3^2 &{}+ T a\varepsilon^3 (x_3^2-x_1^2)x_2 \\ &{} - Tb \varepsilon^3 x_1 x_3^2+ \varepsilon^{3} Tu \end{aligned}\right)\\
%\varepsilon^{3} x_3^{3}
%\end{array}\right)\\&{}= \lambda^{r}(\varepsilon)^3 f(x,u). \nonumber
%\end{align}

% and a terminal cost corresponding to the optimal cost of system \eqref{eq:sys-ex-van-der-pol} linearized around the origin
%and  $\jmath(x):=|x|_P$,  where $P\geq0$ is 
%\begin{equation*}\small
%    P:= \left[\begin{array}{ccc}
%      6.8 & 4    & 0\\
%      4   & 11.5 & 0\\
%      0   & 0    & 0
%    \end{array}\right].
%\end{equation*}
%, and, for similar reasons, we derive that $\jmath$ is  also homogeneous of degree $\mu=2$.
We choose a quadratic stage cost. Hence, for any $x\in\R^{3}$ and $u\in\R$, let $\ell(x,u):=|x|_Q+|u|_R$  with $Q:=\text{diag}(1,1,0)\geq0$ and $R:=1$.
By item (i) of Proposition \ref{lem-sass-ell}, $\ell$ 
is homogeneous of degree $\mu=2$.  Hence, SA\ref{sass-homogeneity-ell} is verified with $\mu=2$. Note that there is no  cost associated to the component $x_3$ in either the stage   or the terminal cost  as we wish for the incurred  cost along solutions of $f$ to coincide with the original system  when $x_3=1$. We also need to check that SA\ref{sass-existence-optimal-value-function} holds. We apply \cite[Theorem 2]{Keerthi-Gilbert-tac85} for this purpose. Indeed, items a) to c) of \cite[Theorem 1]{Keerthi-Gilbert-tac85} hold. Furthermore, item $d_3)$ of \cite[Theorem 2]{Keerthi-Gilbert-tac85} is satisfied as $|u|\leq\ell(x,u)$, hence $u\to\infty$ implies $\ell(x,u)\to\infty$ for any $x\in\R^3$.  The last condition  to check is item e)  of  \cite[Theorem 2]{Keerthi-Gilbert-tac85}, i.e., that for any initial condition there is a sequence of inputs  which makes the cost  $J_{\infty,\gamma_1,\gamma_2}$ finite. For any $x\in\R^3$,  we can construct a linearizing feedback that brings the state to the set $\left\{x \in\R^3 : x_1=x_2=0\right\}$ in at most two steps, while the trivial input $\bm{u}_{\infty}=0$ has cost $J_{\infty,\gamma_1,\gamma_2}((0,0,x_3),\bm{u}_{\infty})=0$  for any $x_3\in\R$. Therefore   an infinite sequence with finite cost exists, and by invoking \cite[Theorem 2]{Keerthi-Gilbert-tac85} we conclude that $V_{d,\gamma_1,\gamma_2}$ exists and SA\ref{sass-existence-optimal-value-function} holds.

All the required conditions hold, we can thus apply Proposition \ref{prop-vi}. We will do so by taking $V_0(x):=\underline{V}_0(x):=\overline{V}_0(x):=|x|_P$ for any $x\in\R^3$, where $P\geq0$ is derived as the LQ optimal cost of  the original system linearized around the origin. That is,
$P:= \left((6.8, 4, 0),(4,11.5,0),(0,0,0)\right)$.
Note that $x\mapsto|x|_P$   is  indeed homogeneous of degree $\mu=2$, in view of Proposition \ref{lem-sass-ell}.
%Now, consider the infinite sequence $\bm{v}(x)=\left\{\dfrac{x_3^2x_1}{T^2}(b-1)-x_2\left(\dfrac{2x_3^2}{T}+a(x_3^2-x_1^2)\right), \dfrac{x_1^+}{T^2}(x_3^6+b)\right.\allowbreak\left.{}+\dfrac{a}{T}(x_3^6-(x_1^+)^2),0,\ldots\right\}$ with $x_1^+\!{:=}\!x_3^2(x_1\!{+}\!T x_2)$. It follows $\phi(1,x,v(x)|_1)=\left(x_1 x_3^2 + T x_2 x_3^2,-\dfrac{x_1 x_3^2 + T x_2 x_3^2}{T}, x_3^{3^k}\right)$ and $\phi(k,x,v(x)|_k)=(0,0,x_3^k)$ for $k\geq2$. Then, for any $d>2$, $\gamma_1\in\Rlo$ and $\gamma_2\in\R$,
%\begin{equation}
%\begin{split}
%\MoveEqLeft J_{d,\gamma_1,\gamma_2}(x,\bm{u}) =  \dst\sum_{k=0}^{d-1}\gamma_1^{\gamma_2^k}\ell\Big(\phi(k,x,\bm{u}),u_k\Big) \\ &{}= \gamma_1^{\gamma_2^0}\ell(\phi(0,x,\bm{u}|_0),u_0)+
%\gamma_1^{\gamma_2^1}\ell(\phi(1,x,\bm{u}|_1),u_1)+0\\ &{}=
%\gamma_1^1\left(x_1^2+x_2^2+(v_1(x))^2\right)\\ &\quad {} +\gamma_1^{\gamma_2}\left(\left(x_1 x_3^2 + T x_2 x_3^2\right)^2\left(1+\frac{1}{T^2}\right)+(v_2(x))^2\right),
%\end{split}
%\end{equation}
%which is finite for any $x\in\R^3$ and $T>0$. Hence, the remaining item  e) of  \cite[Theorem 2]{Keerthi-Gilbert-tac85} holds, and we conclude that $V_{d,\gamma_1,\gamma_2}$ exists and SA\ref{sass-existence-optimal-value-function} holds. Hence, we can apply the results of Section \ref{sect:dp}.
%is a sphere in $\R^3$ of radius $1.5$ that
\vspace{-0.2em}\subsection{Numerical implementation}
We first describe the details of the implementation of the new algorithm in Section \ref{subsect:homVI}. We take the compact set $\mathcal{D}:=\{x\in\R^3,|x|=1.5\}$, which  verifies \eqref{eq:mathcal-D-property}. 

Because we do not know how to solve the minimization steps in \eqref{eq:underline-vi} and \eqref{eq:overline-vi} analytically in general, we proceed numerically. We rewrite the problem in polar coordinates and we consider a  finite difference approximation of $N=501^2$ points equally   distributed in $[-\pi,\pi]\times\left[0,\dfrac{\pi}{2}\right]$, representing  the azimuth and elevation angles of the upper-hemisphere of $\mathcal{D}$, and $M=501$ equally distributed quantized inputs in $[-5,5]$.

To evaluate the accuracy of the proposed algorithm, we need to compute $V_i$ given by standard VI with $V_0(x):=|x|_P$. We consider for this purpose   a finite difference approximation of $N=501^2$ points equally   distributed in $\mathcal{X}:=[-1,1]\times[-1,1]\times\{1\}$  for the state space, and $M=501$ equally distributed quantized inputs in $[-3,3]$ centered at $0$. Note that  the classical VI only produces the value function on $\mathcal{X}$, contrary to $\underline{V_1}(x)$ and $\overline{V_1}(x)$ which are   defined for any $x\in\R^3$. Moreover, this issue is exacerbated at certain points $x\in\mathcal{X}$ \mathieu{where}  $f(x,u)\not\in\mathcal{X}$ for all $u\in[-3,3]$. \mathieu{For standard (approximate) VI schemes to be well defined, an unclear but consequential assessment has to be made to assign values for $V_0(f(x,u))$ when $f(x,u)\not\in\mathcal{X}$.}

\vspace{-0.2em}\subsection{Results}

\begin{figure}
    \centering
    \includegraphics[width=0.86\columnwidth]{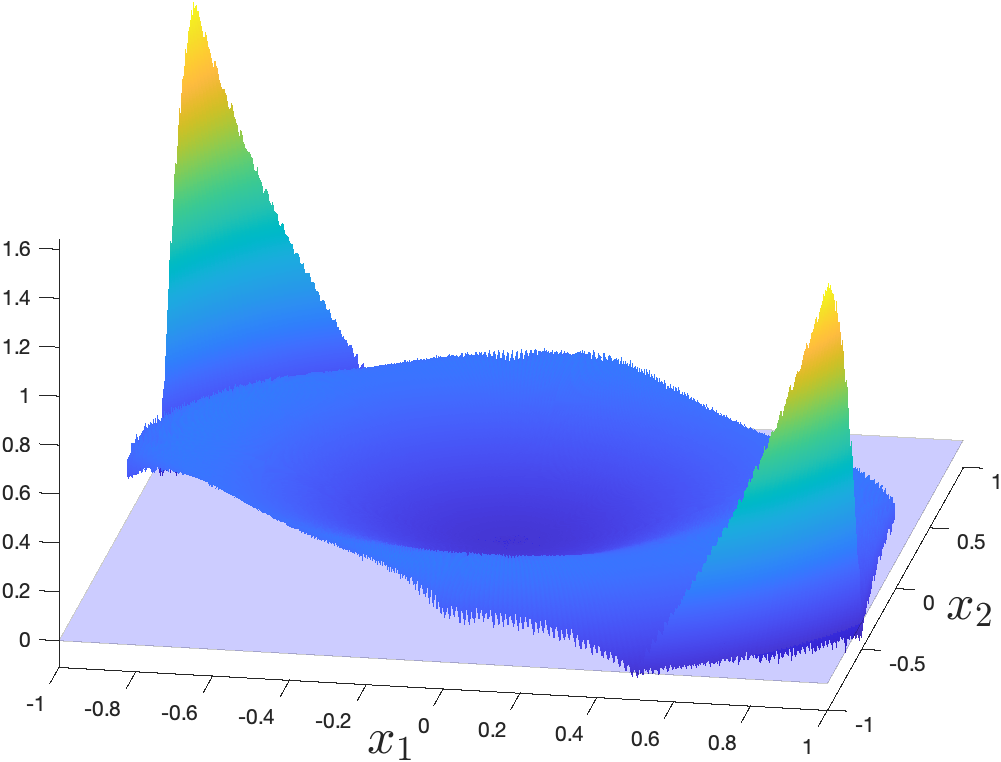}
    \caption{Difference between ${V_1}(x)$ and $\underline{V}_1(x)$ for $x\in\mathcal{X}$}
    \label{fig:under}
    \vspace{-0.25em}
\end{figure}

\begin{figure}
    \centering
    \vspace{-0.25em}
      \includegraphics[width=0.86\columnwidth]{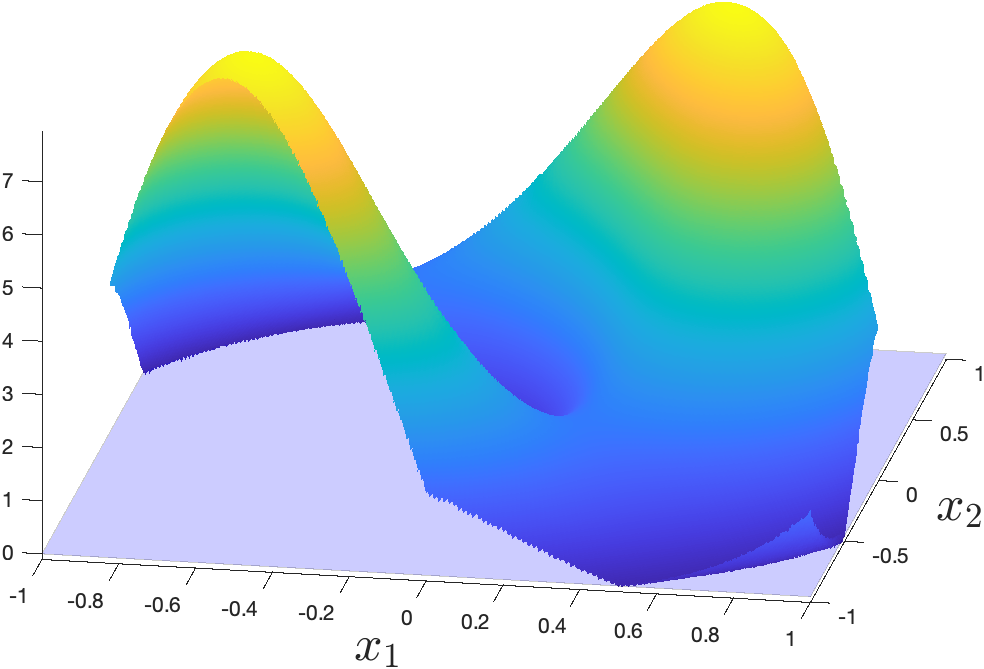}
    \caption{Difference between $\overline{V_1}(x)$ and $V_1(x)$ for $x\in\mathcal{X}$}
    \label{fig:over}
    \vspace{-1em}
\end{figure}

We illustrate the difference between  $V_1(x)$ and $\underline{V_1}(x)$ for $x\in\mathcal{X}$ in \autoref{fig:under}, and similarly   \autoref{fig:over} shows the difference between $\overline{V_1}(x)$  and $V_1(x)$.  As expected, $\underline{V_1}(x)$ and $\overline{V_1}(x)$ under and over estimate $V_1(x)$, respectively. We note from \autoref{fig:under} that the error is small near the intersection of set   $\mathcal{D}$ and $\R^2\times\{1\}$,  which is expected  as it corresponds to points where the minimization in   \eqref{eq:underline-vi} and \eqref{eq:overline-vi} coincides with \eqref{eq:vi}. We have then a circle of radius $\sqrt{1.5^2-1}=1.12$ where the error is small, and  grows when moving to points outside this circle. The error near the origin is also small despite not being part of $\mathcal{D}$, which is explained by the fact that $V_1$,$\underline{V}_1$ and $\overline{V}_1$ vanish at the origin. Similar observations can be made for the upper bound in \autoref{fig:over}, \mathieu{even when} varying the size of the radius  for the compact set $\mathcal{D}$\mathieu{, that is, by taking $\mathcal{D}:=\{x\in\R^3,|x|=r\}$ for different values of $r>0$}.  Hence, by exploiting homogeneity, we have a VI scheme that is calculated only for points in the compact  set $\mathcal{D}$ but  has analytical guarantees  for any state in $\R^{n_x}$, and not just for states in (the forward invariant subset) of $\mathcal{X}$ as in the classical VI implementation.


%\romain{For the cost Mathieu, and so for the satisfaction of SA\ref{sass-homogeneity-ell}-\ref{sass-existence-optimal-value-function}, I let you pick up one that suits you. A quadratic cost would be good, see item (i) of Proposition \ref{prop-sass-f}.}

%\romain{When $\mathcal{D}$ is defined as $\left\{x\in\R^{n_x}\,:\,|x_1|^{\frac{1}{ r_1}}+\ldots+|x_{n_x}|^{\frac{1}{r_{n_x}}}=1\right\}$, given $x\in\R^{n_x}$, we can take $\varepsilon=\dst\left(|x_1|^{\frac{1}{ r_1}}+\ldots+|x_{n_x}|^{\frac{1}{r_{n_x}}}\right)^{-1}$ in that way $\lambda^{r}(\varepsilon)x\in\mathcal{D}$. It would not hurt to double check this.}

\vspace{-0.1em}\section{Concluding remarks}\label{sec:conc}

In this work, we have first adapted the homogeneity notion for discrete-time systems proposed in \cite{Sanchez-et-al-ijrnc19} to be applicable to systems with inputs. We have then provided  several scaling properties of the system, cost and optimal value function. The latter is particularly important as it implies that we only need to solve the original optimal control problem on a given compact manifold of dimension strictly smaller than the dimension of the state space to derive (near-)optimal inputs for any state. We then exploit this observation to derive a new approximation scheme for VI. 

Future work will include further investigation  of the error of approximation when solving VI on the set $\mathcal{D}$ in Section \ref{sect:dp}, as well as other ways to construct upper and lower bounds of the function generated by VI.
\bibliographystyle{plain}
\bibliography{IEEEabrv,main}

%%%%%%%%%%%%%%%%%%%%%%%%%
\end{document}